\documentclass[11pt]{amsart}
\usepackage{latexsym,amsxtra,ifthen,fullpage}

\usepackage{verbatim}
\date{}

\usepackage{amsfonts,amsmath,amsthm}
\usepackage{amssymb,latexsym}
\usepackage[dvips]{epsfig}
\usepackage{amscd}
\usepackage[all]{xy}
\usepackage{enumerate}

\theoremstyle{plain}
\newtheorem*{theorem}{Theorem}
\newtheorem*{lemma}{Lemma}
\newtheorem*{proposition}{Proposition}

\newtheorem*{corollary}{Corollary}
\newtheorem*{cor}{Corollary}

\theoremstyle{definition}

\newtheorem*{remark}{Remark}
\newtheorem*{remarks}{Remarks}

\newcommand{\Zo}{{\mathcal{Z}}}

\newcommand{\cB}{{\mathcal B}}

\newcommand{\cH}{{\mathcal H}}

\newcommand{\cO}{{\mathcal O}}

\newcommand{\cR}{{\mathcal R}}

\newcommand{\cU}{{\mathcal U}}

\newcommand{\B}{{\mathbb B}}

\newcommand{\mg}{\mathfrak{g}}
\newcommand{\mh}{\mathfrak{h}}
\newcommand{\mb}{\mathfrak{b}}

\newcommand{\mC}{\mathbb{C}}

\newcommand{\mZ}{\mathbb{Z}}

\renewcommand{\S}{S(V)}
\newcommand{\T}{S(V^*)}
\newcommand{\ie}{{\it i.e. }}

\newcommand{\la}{\lambda}
\newcommand{\HOM}{\operatorname{Hom}}

\newcommand{\op}{\operatorname}

\newcommand{\ff}{\footnote}

\numberwithin{equation}{section}

\title[Cherednik, Hecke and quantum algebras as Frobenius
extensions]{Cherednik, Hecke and quantum algebras as free
Frobenius and Calabi-Yau extensions}
\author{K.A. Brown, I.G. Gordon, C.H. Stroppel}

\address{Brown: Department of Mathematics,
University of Glasgow, Glasgow G12 8QW, UK}

\email{kab@maths.gla.ac.uk}

\address{Gordon: Department of Mathematics,
University of Glasgow, Glasgow G12 8QW, UK}

\email{ig@maths.gla.ac.uk}

\address{Stroppel: Department of Mathematics,
University of Glasgow, Glasgow G12 8QW, UK}

\email{cs@maths.gla.ac.uk}

\begin{document}
\begin{abstract}
We show how the existence of a PBW-basis and a large enough
central subalgebra can be used to deduce that an algebra is
Frobenius. This is done by considering the examples of rational
Cherednik algebras, Hecke algebras, quantised universal enveloping
algebras, quantum Borels and quantised function algebras. In
particular, we give a positive answer to \cite[Problem
6]{Rouquier} stating that the restricted rational Cherednik
algebra at the value $t=0$ is symmetric.
\end{abstract}
\maketitle

\section{Introduction}

\subsection{}
\label{1.1}
In this note we will consider six types of algebras:
\begin{enumerate}[(I)]
\item the rational Cherednik algebra
  $\op{H}_{0,{\bf c}}$
  associated to the complex reflection group $W$;
  \item the graded (or degenerate) Hecke algebra $\mathbf{H}_{gr}$ associated to
  a complex reflection group $W$;
  \item the extended affine Hecke algebra $\cH$ associated to a finite Weyl
  group $W$;
  \item the quantised enveloping algebra $\cU_\epsilon(\mg)$, at an $\ell$-th root of unity
  $\epsilon$, of a
  semisimple complex Lie algebra $\mg$;
\item the corresponding quantum Borel $\cU_\epsilon(\mg)^{\geq0}$;
\item the corresponding quantised function algebra $\cO_\epsilon[G]$.
\end{enumerate}

These algebras share two important properties: first, they have a regular central subalgebra $\Zo$ over which they
are free of finite rank, second, they - or a closely associated algebra in Case (VI) - have a basis of PBW type.
The purpose of this paper is to show that these two properties are the key tools for defining an associative
non-degenerate $\Zo$-bilinear form for each of these algebras, and hence for deducing Frobenius and Calabi-Yau
properties for the algebras in each class.

\subsection{}
\label{1.2} We prove that each pair $\Zo \subseteq R$ in the
classes (I)-(VI) is a \emph{free Frobenius extension}. The
definition and basic properties are recalled in Section~\ref{frob}
and Section~\ref{Nak} -- in essence, one requires $\HOM_{\Zo}(R,
\Zo) \cong R$ as $(\Zo-R)$-bimodule.
\subsection{}
\label{1.3} When an algebra $R$ is a free Frobenius extension of a
central subalgebra $\Zo$ then $\HOM_{\Zo}(R,\Zo)$ is in fact
isomorphic to $R$ \emph{both} as a left \emph{and} as a right
$R$-module, but not necessarily as a bimodule. However, there is a
$\Zo$-algebra automorphism $\nu$ of $R$, the \emph{Nakayama
automorphism}, such that $\HOM_{\Zo}(R,\Zo) \cong \,{}^1
R^{\nu^{-1}}_{}$ as $R$-bimodules. This automorphism is unique up
to an inner automorphism. We explicitly determine the Nakayama
automorphisms for each case listed above: $\nu$ is trivial
(\ie inner) in cases (I) and (IV); non-trivial in cases (II),
(III) and (V) and (VI).

\subsection{}
\label{1.4} The results summarised in Section~\ref{1.2} have immediate consequences regarding the \emph{Calabi-Yau
property} of the algebras in classes (I) - (VI). The definition and its relevance to Serre duality are recalled in
Section~\ref{2.3}. In particular \cite{CGM}, we get natural examples of so-called Frobenius functors - that is,
functors which have a biadjoint. Frobenius algebras and Frobenius extensions play an important role in many
different areas (see for example \cite{kadison}). They
 give rise to Frobenius functors which are the natural
candidates for constructing interesting topological quantum field
theories in dimension 2 and even 3 (for the latter see for example
\cite{StTQFT}), and also provide connections between
representation theory and knot theory (for example in the spirit
of \cite{KadisonJones}).

\subsection{} \label{1.5}
Let us assume for the moment that $\Zo \subseteq R$ is a free Frobenius
extension with Nakayama automorphism $\nu$. If $I$ is an
ideal of $\Zo$, then it's clear from the definitions that $\Zo/I \subseteq R/IR$ is a free Frobenius extension
with Nakayama automorphism induced by $\nu$. This applies in particular when $I$ is a maximal ideal $\mathfrak{m}$
of $\Zo$; since, for $R$ in classes (I) - (VI), every simple $R$-module is
killed by such an ideal $\mathfrak{m}$, this is relevant to the finite
dimensional representation theory of $R$. Thus $R/\mathfrak{m}R$ is a
Frobenius algebra, which is symmetric provided the automorphism of $R/\mathfrak{m}R$ induced by $\nu$ is inner.

\subsection{}\label{1.6} To define the non-degenerate associative bilinear forms mentioned in Section~\ref{1.1}, we follow in each case
the approach of \cite[Proposition 1.2]{FP} to the study of the inclusion $\Zo \subseteq R$ when $R$ is the
enveloping algebra $U(\mathfrak{g})$ of a finite dimensional restricted Lie algebra $\mathfrak{g}$ over a field
$k$ of characteristic $p > 0,$ and $\Zo$ is the Hopf centre $k \langle x^p - x^{[p]} : x \in \mathfrak{g}
\rangle.$ In the language of the present paper, it is proved there that $\Zo \subseteq U(\mathfrak{g})$ is a free
Frobenius extension, with Nakayama automorphism $\nu$ the winding automorphism of the trace of the adjoint
representation; in particular, $\nu$ is trivial when $U(\mathfrak{g})$ is semisimple. The parallel methods used here
might suggest that an axiomatic approach covering all the cited cases simultaneously might be possible; but we
have not found such a setting.

\subsection{}\label{1.7} The detailed results for classes (I) - (VI) are as follows.
 \begin{enumerate}
 \item (Theorem 3.5 and Corollary 3.6) The rational Cherednik algebra ${\bf H} = \op{H}_{0,c}$ is a
 free Frobenius extension of its central subalgebra $\Zo:=S(V)^W \otimes S(V^*)^W,$ with trivial Nakayama automorphism.
 Consequently ${\bf H}_{\chi}$ is a symmetric algebra for any central character $\chi$ (answering a question of Rouquier,
 \cite[Problem 6]{Rouquier}), and ${\bf H}$ is a Calabi--Yau $\mathcal{Z}-$algebra.
 \item (Theorem~\ref{gradedHecke}) The graded Hecke algebra $\mathbf{H}_{gr}$ associated to a complex
 reflection group $W$ is a free Frobenius
 extension of its centre $\Zo_{gr}:=S(V)^W$, but the Nakayama
 automorphism (which is determined explicitly) is non-trivial.
  \item (Theorem~\ref{affineHecke}) The extended affine Hecke algebra $\cH$ associated to a finite Weyl group $W$ is a free Frobenius
 extension of its centre $\Zo_{\cH}$, but the Nakayama automorphism is non-trivial.
 \item (Theorem~\ref{symmetricU}) The quantised enveloping algebra $\cU_\epsilon(\mg)$ is
 a free Frobenius extension of its $\ell$-centre $\Zo$, with trivial Nakayama automorphism. Consequently,
 $\cU_\epsilon(\mg)_\chi$ is symmetric for any central character $\chi$, and $\cU_\epsilon(\mg)$ is a Calabi-Yau $\Zo -$algebra.
 \item (Theorem~\ref{4.1}) The quantum Borel $\cU_\epsilon(\mg)^{\geq 0}$
 is a free Frobenius extension of its $\ell$-centre $\mathcal{Z}_+$, but the Nakayama
 automorphism (which is determined explicitly) is non-trivial.
 \item (Theorem~\ref{Frobenius}) There is an element $z$ of the central subalgebra $\cO[G]$ of the quantised function algebra $\cO_{\epsilon}[G]$ such that $\cO_{\epsilon}[G][z^{-1}]$ is a free Frobenius extension of $\cO[G][z^{-1}]$ with non-trivial Nakayama automorphism. The open set $\mathcal{O}_z = \{ g\in G: z\notin \mathfrak{m}_g\}$ meets every torus orbit of symplectic leaves in $G$. Thus, for any $g\in G$, the algebra $\cO_\epsilon[G]/\mathfrak{m}_g \cO_\epsilon[G]$ is Frobenius but not, in general, symmetric.
 \end{enumerate}

 \subsection{}\label{1.8} There is some overlap between this paper and \cite{Braun}, a preliminary version of which we
 received while this paper was being written. The methods used in the two papers are completely different, and
 indeed complementary.

\subsection{}\label{1.9} In the following rings are always assumed to be
 unitary and, if not stated otherwise, modules are {\it left} modules. For any
 ring $S$ we denote by $\HOM_S(-,-)$, $\HOM_{-S}(-,-)$ and $\HOM_{S-S}(-,-)$
 the morphism spaces in the category of (left) $S$-modules, right $S$-modules
 and $S$-bimodules respectively. Our algebras are all over $\mC$; undoubtedly this
 hypothesis could be weakened. We abbreviate $\otimes=\otimes_\mC$.

\section{Frobenius and Calabi-Yau extensions}
\label{Frobeniusext}
\subsection{Definition}\label{frob}
We first recall some basics on Frobenius extensions.  For more details we refer for example to \cite{BF},
\cite{K}, \cite{NT}, \cite{Pareigis}. A ring $R$ is a {\it free Frobenius extension} (of the first kind) over a
subring $S$, if $R$ is a free $S$-module of finite rank, and there is an isomorphism of $R-S$-bimodules
$F:R\longrightarrow \HOM_{S}(R,S)$. (The bimodule structure on the latter is defined as $r.f.s(x)=f(xr)s$ for $r,
x\in R$, $s\in S$, $f\in \HOM_{S}(R,S)$.) Equivalently, $R$ is a free right $S$-module of finite rank, and there
is an isomorphism of $S-R$-bimodules $G:R\longrightarrow \HOM_{-S}(R,S)$ (\cite[Proposition 1]{NT}). The existence
of $F$ provides a non-degenerate associative $S$-bilinear form $\B :R\times R\rightarrow S$, defined by $\B
(r,t)=F(t)(r)$ for all $r,t \in R$. Given a basis $r_i$, $1\leq i\leq n$ of $R$ as an $S$-module, we find elements
$r^i$, $1\leq i\leq n$ such that $\B (r_i, r^j)=\delta_{i,j}$ because $F$ is surjective. The two ordered sets
$\{r_i: 1\leq i\leq n\}$ and $\{r^i: 1\leq i\leq n\}$ form a dual free pair (in the sense of \cite[Section
1]{BF}). Conversely, the existence of a non-degenerate associative bilinear form $\B :R\times R\rightarrow S$
together with a dual free pair implies that $R$ is a free Frobenius extension of $S$ with defining isomorphism $F$
given by $F(t)(r)=\B(r,t)$ (see \cite[Section 1]{BF}).

\subsection{The Nakayama automorphism}
\label{Nak} We recall some ideas from \cite{K}. Suppose for the
rest of this section that $R$ is a free Frobenius extension of $
\mathcal{Z}$, with $ \mathcal{Z}$ now contained in the centre of
$R$. The isomorphisms $F$ and $G$ defined in 2.1 induce
isomorphisms of left respectively right $R$-modules
\begin{eqnarray}
\label{eq:FG}
  \begin{array}[thl]{ccccccc}
  R&\cong&\HOM_{ \mathcal{Z}}(R, \mathcal{Z})&=&\HOM_{ \mathcal{Z}- \mathcal{Z}}(R, \mathcal{Z})&=&RF(1)\\
  R&\cong&\HOM_{- \mathcal{Z}}(R, \mathcal{Z})&=&\HOM_{ \mathcal{Z}- \mathcal{Z}}(R, \mathcal{Z})&=&G(1)R.
  \end{array}
\end{eqnarray}
One can show \cite[Section 2 (4)]{K} that $h:=F(1)=G(1)$ as elements of $\HOM_{ \mathcal{Z}- \mathcal{Z}}(R,
\mathcal{Z})$. Thus we get a well-defined $\mathcal{Z}$-algebra automorphism $\nu:R \longrightarrow R$, defined by
$rh = h\nu (r)$ for all $r\in R$. An easy calculation shows that
\begin{eqnarray*}
 \B(x,y)&=&\B(\nu(y), x)
\end{eqnarray*}
for $x$, $y\in\cB$. The automorphism $\nu$ is called the {\it Nakayama
  automorphism} (with respect to $F$, $\B$, or $G$). It's clear that $\nu$ is uniquely determined up to an inner automorphism of $R$ by the
pair
  $\mathcal{Z} \subseteq R$. It therefore makes sense to speak about the Nakayama
  automorphism attached to a free Frobenius extension. We call the extension
  {\it symmetric} if the Nakayama automorphism is inner.

Thanks to our assumption on $\mathcal{Z}$, there is now also a \emph{right} $R$-action on
  $\HOM_{\mathcal{Z}}(R,\mathcal{Z})$, given by $fr(-) = f(r-)$ for $r \in R$ and $f \in \HOM_{\mathcal{Z}}(R,\mathcal{Z})$.
   Let $_{}^{1}R_{}^{\nu^{-1}}$
  be the ring $R$ considered as an $R$-bimodule, but with its right $R$-module
  structure twisted by $\nu^{-1}$. Then the $R-\mathcal{Z}$-bimodule isomorphism $F$ is in fact an
  isomorphism of $R$-bimodules
  \begin{eqnarray}
    \label{eq:CY}
    _{}^{1}R_{}^{\nu^{-1}}\cong\HOM_{\mathcal{Z}}(R,\mathcal{Z}),
  \end{eqnarray}
since $F(r\nu^{-1}(x))(y)=F(\nu^{-1}(x))(yr)=\B(yr,\nu^{-1}(x))=\B(x,yr)$ and
$(F(r)x)(y)=F(r)(xy)=F(yr)(x)=\B(x,yr)$ for all $x,y,r \in R$.

\subsection{} \label{hyp} We now highlight a condition which will allow us to prove
that algebras are free Frobenius extensions. For this we let $R$
be free  with a finite basis $\mathcal{B}$ over an affine central
subalgebra $\mathcal{Z}$. The condition is:

\medskip
\noindent {\bf Hypothesis:} There exists a $\mathcal{Z}$--linear functional $\Phi : R\rightarrow \mathcal{Z}$ such
that for any non-zero $a = \sum_{b \in \mathcal{B}} z_b b\in R$ there exists $x\in R$ with $\Phi(xa) = uz_b$ for
some unit $u\in \mathcal{Z}$ and some non-zero $z_b\in \mathcal{Z}$.

\begin{proposition} \label{cruc} Let $R$ be a finitely generated free $\mathcal{Z}$--module with a basis $\mathcal{B}$ which satisfies the above hypothesis. Then $R$ is a free Frobenius extension of $\mathcal{Z}$ and for any maximal ideal $\mathfrak{m}$ of $\mathcal{Z}$, the finite dimensional quotient $R/\mathfrak{m}R$ is a finite dimensional Frobenius algebra.
\end{proposition}
\begin{proof}
Let $\theta : R \rightarrow \HOM_\mathcal{Z}(R,\mathcal{Z})$ be the $R-\mathcal{Z}$-bimodule homomorphism defined
by $\theta(a)(a') = \Phi (a'a)$. Clearly $\theta$ is an injection since if $a\in R$ is non-zero then the
displayed hypothesis implies that $\theta(a)(x) \neq 0$. Thus we have a short exact sequence
\begin{equation} \label{frobsequence} 0\rightarrow R \rightarrow \HOM_\mathcal{Z}(R,\mathcal{Z}) \rightarrow C
\rightarrow 0\end{equation} of $R-\mathcal{Z}$-bimodules, where $C$ is the cokernel of $\theta$. We will prove
that $C=0$ after showing that $\theta$ induces a Frobenius structure on each finite dimensional quotient
$R/\mathfrak{m}R$.

Fix an arbitrary maximal ideal $\mathfrak{m}$ of $\mathcal{Z}$ and consider the mapping $$\overline{\theta} : \frac{R}{\mathfrak{m}R} \longrightarrow \HOM_\mathcal{Z}(R,\mathcal{Z}) \otimes_\mathcal{Z} \frac{\mathcal{Z}}{\mathfrak{m}}$$ which sends $a+\mathfrak{m}R$ to $\theta(a) \otimes 1$. Let $$\iota : \HOM_\mathcal{Z}(R,\mathcal{Z}) \otimes_\mathcal{Z} \frac{\mathcal{Z}}{\mathfrak{m}} \longrightarrow \HOM_{\mC}( \frac{R}{\mathfrak{m}R}, \mC)$$ be the isomorphism sending $\psi\otimes 1$ to the mapping $(a + \mathfrak{m}R \mapsto \psi(a) + \mathfrak{m})$.

We claim that composition $\iota \overline{\theta}$ is an isomorphism. To prove this, we will show that $\iota \overline{\theta}$ is injective; then, since both the domain and codomain are vector spaces of the same dimension, the claim will follow. By construction, $$\iota\overline{\theta} ( a+ \mathfrak{m}R ) (a' + \mathfrak{m}R) = \Phi (a' a) + \mathfrak{m}.$$ Therefore, if $a+\mathfrak{m}R \in \ker \iota\overline{\theta}$ then $\Phi (a'a) \in \mathfrak{m}$ for all $a'\in R$. We assume that $a\neq 0$. Then, by hypothesis, if we write $a = \sum z_b b$, we can find $x\in R$ such that $\Phi (xa) = u z_b$ for some unit $u$ and some non-zero $z_b$. Thus $z_b \in \mathfrak{m}$. Now $a$ and $a-z_bb$ have the same image in $R/\mathfrak{m}R$ so we can replace $a$ by $a-z_bb$. Repeating this procedure shows that $a \in \mathfrak{m}R$ and hence that $\iota\overline{\theta}$ is injective.

As a first consequence we see that $\iota\overline{\theta}$ induces an $R/\mathfrak{m}R$-isomorphism $R/\mathfrak{m}R \cong (R/\mathfrak{m}R)^*$ so $R/\mathfrak{m}R$ is Frobenius. We also deduce that $\overline{\theta}$ is an isomorphism, and so from \eqref{frobsequence} we see $C\otimes_\mathcal{Z} \mathcal{Z}/\mathfrak{m}\mathcal{Z} = 0$. Since this is true for an arbitrary maximal $\mathfrak{m}$ of $\mathcal{Z}$ and $C$ is finitely generated over $\mathcal{Z}$, it follows that $C = 0$. Hence $\theta : R \longrightarrow \HOM_\mathcal{Z}(R,\mathcal{Z})$ is an isomorphism and so $R$ is a free Frobenius extension of $\mathcal{Z}$.
\end{proof}
\subsection{Calabi-Yau algebras}
\label{2.3} Let $d$ and $n$ be non-negative integers and let $R$ be a ring which has a commutative noetherian
central subring $C$ of Krull dimension $d$, over which $R$ is a finitely generated module. Following for example
\cite{IR}, we say that $R$ is a \emph{Calabi-Yau $C$-algebra of dimension $n$} if, for all $X,Y \in \mathcal{D}^b
(\mathrm{Mod} (\mathrm{fl-}R))$, the bounded derived category of $R-$modules of finite length, there is a natural isomorphism
\begin{eqnarray*}
\HOM_{\mathcal{D}(\mathrm{Mod} (R))}(X,Y[n]) \cong D\HOM_{\mathcal{D}(\mathrm{Mod} (R))}(Y,X). \end{eqnarray*}
Here, $D$ denotes the \emph{Matlis duality} functor $D = \HOM_C (-,E),$ where $E$ is the direct sum of the
$C-$injective hulls of the simple $C-$modules. The following proposition is an immediate consequence of
\cite[Theorems 3.1 and 3.2]{IR}, once we note that if $C$ is regular then the Cohen-Macaulay $C-$modules coincide
with the projective $C-$modules.

\begin{proposition}
Let $C$, $R,$ $n$ and $d$ be as above, and suppose that $C$ is a regular domain. Then $R$ is a Calabi-Yau
$C-$algebra of dimension $n$ if and only if $n=d,$ $R$ has finite global dimension, $R$ is a projective
$C-$module, and $\HOM_C(R,C)$ is isomorphic to $R$ as $R-R-$bimodules. In this case, $R$ has global dimension $d.$
\end{proposition}

\subsection{Hopf algebras}1. When $H$ is a Hopf algebra which is a finite module over a central affine Hopf
subalgebra $\mathcal{Z},$ Hopf-algebraic methods can be used to
deduce that $H$ is a Frobenius extension of $\mathcal{Z}$. The
result is due to Kreimer and Takeuchi \cite[Theorem 1.7]{KT}; the
arguments are sketched in \cite[Section III.4]{BGo}. This provides
an alternative approach to the algebras in classes (IV), (V) and
(VI), but this does not provide an explicit description of the
bilinear form, nor does it give immediate access to the Nakayama
automorphism.

2. The concept of the Nakayama automorphism was introduced also in a recent paper on noetherian Hopf algebras by
Brown and Zhang \cite{BZ}. They showed that many noetherian Hopf algebras $H$ (including all those which are
finite modules over their centres) have a rigid dualizing complex $R$ which is isomorphic (in the derived category
of bounded complexes of $H-$bimodules) to ${}^{\hat{\nu}}H^1 [d]$; here, $d$ is the injective dimension of $H$,
$[d]$ denotes the shift, and $\hat{\nu}$ is a certain algebra automorphism of $H$ which Brown and Zhang called the
Nakayama automorphism. The automorphism $\hat{\nu}$ is trivial on the centre of $H$ and is uniquely determined by
$H$, up to an inner automorphism.

When both usages of the term "Nakayama automorphism" are in play, they define the same map (bearing in mind that
both definitions are only unique up to an inner automorphism of the algebra). To see this, suppose that $H$ is a
free Frobenius extension of a smooth affine central subalgebra $\mathcal{Z}$, (as is the case for the algebras of
(IV), (V) and (VI)). Then the injective dimension $d$ of $H$ equals the Krull dimension (of $H$ and of
$\mathcal{Z}$). Thus the rigid dualizing complex of $\mathcal{Z}$ is $\mathcal{Z}[d]$, and, by \cite[Proposition
5.9]{Ye}, \cite[Example 3.11]{YZ}, $H$ has rigid dualizing complex
$\mathrm{RHom}_{\mathcal{Z}}(H,\mathcal{Z}[d])$. From the free Frobenius property of $H$, and (2.2), we deduce
that this latter complex is isomorphic to ${}^{\nu}H^1[d],$ where $\nu$ denotes the Nakayama automorphism of the
present paper. By the uniqueness of the rigid dualizing complex of $H$ \cite[Proposition 8.2]{VdB}, it follows
that $\hat{\nu} = \nu$ up to an inner automorphism, as claimed.

\section{The rational Cherednik algebra}
\label{6}
In this section we show that the rational Cherednik algebra
${\bf H}$ is a Frobenius
extension of its (what we call) bi-invariant centre, with trivial Nakayama
automorphism, so that the reduced Cherednik algebras
${\bf H}_\chi$ are symmetric.
\subsection{Rational Cherednik Algebras}
\label{6.1}
Let $W$ denote an irreducible complex reflection group with identity element $e$ and set of complex reflections $S$. We fix $V$, a complex reflection representation of $W$, and set $n = \dim V$. Let $c$ be
a conjugation invariant complex function on $S$. For $s\in
S$ let $\alpha_s$ (respectively $\check{\alpha}_s$) be a linear functional
on $V$ (respectively $V^*$) which vanishes on the reflection hyperplane
for $s$; we normalise these by the condition $\langle \alpha_s,\check{\alpha}_s\rangle=2$. {\it The
rational Cherednik algebra $\mathbf{H}=\mathbf{H}_{0,c}$} is the $\mC$-algebra
generated by $\{w\in W, x\in V, y\in V^*\},$ with defining
relations
\begin{eqnarray}
  &wxw^{-1}={}^{w}x,\quad wyw^{-1}={}^wy,&\label{def1}\\
&[x,x']=0,\quad [y,y']=0,&\label{def1.5}\\
  &\left[x,y\right]=\sum_{s\in S}c(s)\langle
  y, \alpha_s\rangle\langle\check{\alpha}_s,x\rangle s&,\label{def2}
\end{eqnarray}
for $x, x'\in V$, $y, y'\in V^*$ and $w \in W.$ These are the
algebras $\mathbf{H}_{0,c}$ from \cite[p.251]{EG}.

\subsection{The PBW-basis}
The algebra $\mathbf{H}$ has a PBW-property in the following sense: multiplication induces an isomorphism $$S(V)\otimes_\mC\mC W\otimes_\mC
S(V^*)\tilde{\longrightarrow}\mathbf{H}$$ of vector spaces (see~\cite[Theorem 1.3]{EG}). In
particular, there is a PBW-basis given by the elements of the set
$\cB_{\mathbf{H}}=\{fwg\}$, where $w\in W,$ $f$ runs through a
homogeneous basis of $\S,$ and $g$ runs through a homogeneous
basis of $\T$.

For $f$ in $S(V)$ or $\T$ we write $|f|$ for the degree of $f$.
For $i \in \mZ_{\geq 0}$ let $\cB_{<i}$ be the span of all
PBW-basis elements of the form $fxg$, where $f\in S(V)$, $x\in W$
and $g\in\T,$ such that $f$ and $g$ are homogeneous with $|f| + |g| < i$: this induces a filtration of $\mathbf{H}$. Moreover, the commutation relation \eqref{def2} shows that \begin{equation} \label{commforus} \text{$|[f,g]| \leq |f| + |g| - 2$ for all homogeneous $f\in S(V)$ and $g\in \T$.}\end{equation}

\subsection{The central subalgebra}
\label{6.2} The algebra $\mathbf{H}=\mathbf{H}_{0,c}$ has a large centre $Z({\bf H})$, isomorphic to the so-called
spherical subalgebra (\cite[Theorem 3.1, Theorem 7.2]{EG}). In particular, $Z({\bf H})$ contains the {\it bi-invariant centre} $$\Zo \quad :\quad \S^W\otimes\T^W.$$ Now $\S$ (respectively $\T$) is a free $\S^W$-module (respectively
$\T^W$-module) of rank $|W|$, see \cite[V.18.3]{Kane} for example. A basis can be obtained by taking arbitrary homogeneous preimages of any homogeneous
basis of the coinvariant algebra $A:=\S/(\S^W_+)$.

Then $A$ is a local Frobenius algebra thanks to \cite[Proposition
VII.26.7]{Kane} and its associated bilinear form is easy to
describe. To do this, recall that the homogeneous component $A_N$
of $A$ of highest degree has dimension one and is skew invariant
for the action of $W$ on $V$, \cite[20.3, Propositions A and
B]{Kane}. Let $\pi : A \rightarrow A_N$ be the projection map with
$\pi (A_i) = 0$ for $i \neq N.$ Then the bilinear form is given by
$$B(\overline{a}, \overline{a'}) = \pi(\overline{a'}\overline{a}).$$ Similar statements apply to $\T/(\T^W_+)$: it is Frobenius and its
highest degree component is skew invariant for the action of $W$ on $V^*$. Below, we shall use the notation
$\epsilon_V, \epsilon_{V^*}$ for these two one-dimensional representations of $W$.

\subsection{}
\label{dualbases}
We fix a pair of homogeneous dual bases $\{\overline{{\bf a}}_i:1\leq i\leq
    |W|\}$, $\{{\overline{\bf a}}^i:1\leq i\leq |W|\}$ for $\S/(\S^W_+)$, and  a
    pair of homogeneous dual bases $\{{\overline{\bf
    b}}_i:1\leq i\leq |W|\}$, $\{{\overline{\bf b}}^i:1\leq i\leq |W|\}$ for $\T/(\T^W_+)$. Then we lift them to homogeneous $S(V)^W$-bases, $\{{\bf
a}_i:1\leq i\leq |W|\}$, $\{{\bf
  a}^i:1\leq i\leq |W|\}$ of $\S$, and homogeneous $\T^W$-bases $\{{\bf
  b}_i:1\leq i\leq |W|\}$,
$\{{\bf b}^i:1\leq i\leq |W|\}$ of $\T$. We set $|{\bf a}_i| d_i$ and $|{\bf b}_i|= e_i$; then $|{\bf a}^i| = N-d_i$ and $|{\bf
b}^i| = N - e_i$. Let ${\bf
  a}_{max}$ and ${\bf b}_{max}$ be the elements of
maximal degree $N$ amongst the ${\bf a}_i$ and ${\bf b}_i$ respectively.

\subsection{The functional}
For $f\in\S$ let ${\bf a}_{max}(f)$ be the coefficient of ${\bf
a}_{max}$ when $f$ is expressed in the chosen $S(V)^W$-basis of
$S(V)$. Similarly, we define ${\bf b}_{max}(g)$ for $g\in\T$.
Thanks to the PBW-property, $\mathbf{H}$ is a free $\Zo$-module of
finite rank with basis  $$\cB_{\mathbf{H}} \quad := \quad \{{\bf
a}_{i}w{\bf b}_{j} : w\in W, 1\leq i,j\leq |W|\}.$$
We define a
$\Zo$-linear map
\begin{eqnarray*}
  \Phi:\quad\quad\quad\mathbf{H}&\longrightarrow&\Zo\\
  \cB_{\mathbf{H}}\ni {\bf a}_iw{\bf b}_j&\longmapsto&
  \begin{cases}
1 \quad &\text{if $i=j=max$ and $w=e$},\\
0 \quad &\text{otherwise}.
  \end{cases}
\end{eqnarray*}

\label{6.3}
\begin{lemma}
\label{main} The functional $\Phi$ above satisfies Hypothesis~\ref{hyp}.
\end{lemma}

\begin{proof}
Let $a = \sum_{b\in \mathcal{B}_{\mathbf{H}}} z_b b$ be a non-zero
element of $\mathbf{H}$.
 Pick $b = {\bf a_i} w {\bf b_j} \in \mathcal{B}_{\mathbf{H}}$ of maximal degree $|{\bf a_i}|+|{\bf b_j}|$
 such that $z_b \neq 0,$ and set $x = {\bf b}^jw^{-1} {\bf a}^i$. We claim that this choice of $x$ satisfies Hypothesis~\ref{hyp}.

For indices $i', j'$ and for $u\in W$ we have, by \eqref{def1} and \eqref{commforus}, \begin{eqnarray*} {\bf
a}_{i'} u{\bf b}_{j'}x={\bf a}_{i'} u{\bf b}_{j'} {\bf b}^jw^{-1}{\bf a}^i &=& {\bf a}_{i'}\cdot uw^{-1}\cdot
{}^w( {\bf b}_{j'}{\bf b}^j) {\bf a}^i \\ & = & {\bf a}_{i'} \, ({}^{uw^{-1}}{\bf a}^i) \cdot uw^{-1} \cdot
{}^{w}({\bf b}_{j'}{\bf b}^j) + \text{lower order terms}. \end{eqnarray*} Since $b$ was chosen to have maximal
degree it follows that if ${\bf a}_{i'} u{\bf b}_{j'}$ appears in the expansion of $a$, then the lower order terms
in the above expression have total degree less than $d_i + e_j + (N-d_i) + (N-e_i) = 2N$. Therefore we find that
\begin{eqnarray}\nonumber
\Phi ({\bf a}_{i'} u{\bf b}_{j'}x)\Phi ({\bf a}_{i'} u{\bf b}_{j'} {\bf
b}^jw^{-1}{\bf a}^i ) &=& \Phi (  {\bf a}_{i'} \, ({}^{uw^{-1}}{\bf a}^i) \cdot uw^{-1} \cdot {}^{w}({\bf
b}_{j'}{\bf b}^j))\\ \label{gut} & = & \delta_{u,w}\Phi ( {\bf a}_{i'} {\bf a}^i \cdot {}^{w}({\bf b}_{j'}{\bf
b}^j)).\end{eqnarray} By definition of the dual basis we have, for $i,i',j,j' = 1, \ldots , N,$
\begin{eqnarray*}{\bf a}_{i'}{\bf a}^{i} = (\delta_{i,i'}+ r_{max}){\bf a}_{max} + \sum_{k\neq max} r_k {\bf a}_k \quad \mathrm{
and } \quad {\bf b}_{j'}{\bf b}^{j} = (\delta_{j,j'}+ r'_{max}){\bf b}_{max} + \sum_{k\neq max} r'_k {\bf
b}_k\end{eqnarray*} for some $r_{max}, r_k \in S(V)^W$ and $r'_{max}, r'_k \in S(V^*)^W$. Consideration of
polynomial degrees in the above expressions shows that $r_k \in (S(V)^W_+)$ for $k = max,$ and for all $k$ when $i
= i',$ and that $r'_k \in (S(V^*)^W_+)$ for $k = max,$ and for all $k$ when $j = j'.$ Substituting in (\ref{gut})
we find that there exists $0\not=c\in\mC$ such that
\begin{eqnarray}
\label{key} & \Phi ( {\bf a}_{i'} {\bf a}^i \cdot {}^{w}({\bf b}_{j'}{\bf b}^j)) = \\ \nonumber
&c\Phi((\delta_{i,i'}+ r_{max}){\bf a}_{max} + \sum_{k\neq max} r_k {\bf a}_k)((\delta_{j,j'}+ r'_{max}){\bf
b}_{max} + \sum_{k\neq max} r''_k {\bf b}_k)),
\end{eqnarray} where $r''_k \in S(V^*)^W$ and $r''_k \in (S(V^*)_+^W)$ when $j = j'.$ We claim that (\ref{key}) is 0 except when $(i',j') = (i,j).$ To see this, suppose that $(i',j')$ is not equal to
$(i,j)$, but (\ref{key}) is non-zero. Our choice of $b$ to have maximal degree with $z_b \neq 0$ forces
\begin{eqnarray}
\label{count} d_{i'} + e_{j'} = d_i + e_j, \end{eqnarray} since otherwise the degree of ${\bf a}_{i'} {\bf a}^i
\cdot {}^{w}({\bf b}_{j'}{\bf b}^j)$ is strictly less than $2N$, and hence can't involve ${\bf a}_{max}{\bf
b}_{max}.$

Suppose first that $i' \neq i$ and $j' \neq j.$ Then (\ref{key}) becomes
\begin{eqnarray}
\label{dubya}  \Phi ( {\bf a}_{i'} {\bf a}^i \cdot {}^{w}({\bf b}_{j'}{\bf b}^j)) = r_{max}r'_{max} \Phi({\bf
a}_{max}{\bf b}_{max}).
\end{eqnarray}
But $r_{max}, r'_{max}$ are in the ideals of positive degree invariants, and so have strictly positive degrees if
they are not 0. Thus, comparing degrees in (\ref{dubya}), using (\ref{count}), shows that (\ref{dubya}) is 0 in
this case.

Suppose now that $i = i'$ but that $j \neq j'.$ Then, by (\ref{count}), $e_{j'} = e_j.$ Therefore
\begin{eqnarray}\label{nixon} {\bf b}_{j'}{\bf b}^{j} =  r'_{max}{\bf b}_{max} + \sum_{k\neq max} r'_k {\bf b}_k,
\end{eqnarray}
and in this equation $r'_{max}=0$, since otherwise it has strictly positive degree, contradicting the
homogeneity of degree $N$ of (\ref{nixon}). Hence (\ref{key}) becomes

$$ \Phi ( {\bf a}_{i} {\bf a}^i \cdot {}^{w}({\bf b}_{j'}{\bf b}^j)) \Phi(({\bf a}_{max} + \sum_{k\neq max} r_k {\bf a}_k)(\sum_{k\neq max} r''_k {\bf b}_k)) = 0.$$ A similar argument
applies if $i' \neq i$ but $j' = j.$ Thus the claim is proved. Therefore $$ \Phi({\bf a}_{i'} u{\bf b}_{j'} {\bf
b}^jw^{-1}{\bf a}^i ) = \delta_{u,w}\delta_{i,i'}\delta_{j,j'}\epsilon_{V^*}(w). $$ It follows that, with $x {\bf b}^j w^{-1} {\bf a}^i,$
$$ \Phi (ax) = z_b \epsilon_{V^*}(w) $$ where $b = {\bf a}_i w {\bf b}_j,$ confirming Hypothesis 2.3.
\end{proof}

\subsection{The theorem for Cherednik algebras} Define the form $\B$ for $\mathbf{H}$ by $\B(a,b) = \Phi (ab),$ for $a,b \in \mathbf{H}.$ We can now deduce the

\begin{theorem}
  The rational Cherednik algebra $\mathbf{H}$ is a symmetric
  Frobenius extension of its central subalgebra $\Zo = S(V)^W \otimes S(V^*)^W.$
\end{theorem}

\begin{proof} It is immediate from Lemmas \ref{6.3} and 2.3 that ${\bf H}$ is a free Frobenius
extension of $\mathcal{Z}$ with form $\B$ as defined above.
Therefore it remains only to prove that the Nakayama automorphism
for ${\bf H}$ is inner.

We verify that $\B(Y,x)=\B(x,Y)$, where $Y\in\cB_{\mathbf{H}}$ and $x\in W$ or $V$ or $V^*$, since $W$, $V$ and $V^*$ generate $\mathbf{H}$ as a $\Zo$-algebra. Let $fwg$ be a typical element from
  $\cB_{\mathbf{H}}$. First, let $x\in W.$ Then
  \begin{eqnarray}
    \B(fwg,x)=\Phi(fwgx)&=&\Phi(f \cdot wx\cdot {}^{x^{-1}}g)\label{1}\\
    &=&\epsilon_{V^*}(x^{-1})\Phi(f\cdot wx \cdot g)\label{2}\\
    &=&\epsilon_{V^*}(x^{-1})\Phi(f\cdot xw\cdot g)\label{3}\\
    &=&\epsilon_{V}(x^{-1})\epsilon_{V^*}(x^{-1})\Phi({}^xf \cdot xw \cdot g)\label{4}\\
    &=&\Phi(xfwg)\label{5} = \B(x,fwg)
  \end{eqnarray}
The equalities~\eqref{1} follow from the definition of $\B$ and
the defining relations~\eqref{def1} of $\mathbf{H}$. To see the
formulas~\eqref{2} and~\eqref{4} note that $x({\bf
a}_{max})=\epsilon_V(x){\bf a}_{max}+h$, where $h \in S(V)$ with
${\bf a}_{max}(h)=0$. Similarly for ${\bf b}_{max}$, and then
invoke the definition of $\Phi$. The equality~\eqref{3} is true
because both sides of the equation are trivial unless $x=w^{-1}$,
in which case we have $xw=wx$. The relation~\eqref{5} holds
because of the defining relations of ${\mathbf H}$ and thanks to
the fact that $\epsilon_V(x)=\epsilon_{V^*}(x)^{-1}$. Finally, the last equation is clear by definition of $\B$, and
hence $\B(fwg,x)= \B(x,fwg)$ holds.

If $a\in V$ we get
\begin{eqnarray}
    \B(fwg,a)&=&\Phi(fwga) \nonumber\\ & =& \Phi(fwag)\label{11}\\
    &=&\Phi(f \, {}^wa \, wg)\nonumber\\
    &=&\Phi(fawg)\label{22}\\
    &=&\Phi(afwg)=\B(a,fwg).\nonumber
  \end{eqnarray}
The equality in \eqref{11} arises since the degree of $fwga$ and
$fwag$ is $|f| + |g| + 1$ and so both sides are zero unless $|f| +
|g| \geq 2N-1$. In the case $|g| = N$ or $N-1$ then $|[a,g]| < N$
by \eqref{commforus}. This then means that $\Phi(fwga) = \Phi(
fwag - fw[a,g]) = \Phi ( fwag)$, as required. The
equality~\eqref{22} is true, because we have zero on both sides if
$w\not=e$. Hence $\B(fwg,a)=\B(a,fwg)$ holds. If $b\in V^*$ the
argument is similar, so we leave it to the reader.

Therefore we get $\B(x,y)=\B(y,x)$ for any $x, y\in\mathbf{H}$, which means $\B$ is symmetric.
\end{proof}

\subsection{Consequences}
\label{conseq} Given a maximal ideal $\mathfrak{m}_{\chi}$ of $\mathcal{Z}$ we define the {\it reduced Cherednik
algebra} to be the $|W|^3$-dimensional algebra $${\bf H}_{\chi} \equiv \frac{\bf H}{\mathfrak{m}_{\chi} \bf H}.$$
Thanks to \cite{IGG} these algebras control a great deal of the geometry associated to the centre of {\bf H}. The
following corollary is immediate from Theorem 3.5 and the discussion in 2.4, after we have noted that ${\bf H}$
has finite global dimension by \cite[page 276]{EG}. The first part (for the case when $\mathfrak{m}_{\chi}$ is
$(S(V)^W \otimes S(V^*)^W)_+$) answers \cite[Problem 6]{Rouquier}.

\begin{corollary} \begin{enumerate}
\item The reduced Cherednik algebras ${\bf H}_{\chi}$ are
symmetric, with dual bases the images of the bases $\mathcal{B} = \{{\bf a_i}w{\bf b_j}\}$ and $\mathcal{B}' \{{\bf a^i}w{\bf b^j}\}$ defined in Section~\ref{6.2} and \ref{dualbases}. \item ${\bf H}$ is a Calabi--Yau  $\Zo-$algebra of dimension $2
\dim (V).$
\end{enumerate}
\end{corollary}

\section{The graded Hecke algebra}
In this section we show that the graded Hecke algebra
${\bf H}_{gr}$ is a Frobenius
extension of its invariant centre, with non-trivial Nakayama automorphism,
so that the reduced graded Hecke algebras
${{\bf H}_{gr}}_\chi$ are Frobenius but not, in general, symmetric.

\subsection{Graded Hecke algebras}
\label{GradedHecke}
As in the previous section let $W$ be an irreducible complex reflection group
with identity $e$, and $V$ the defining complex reflection representation of $W$. Let $\mathbf{H}_{gr}$ be the associative algebra
generated by $V$ and $\mC W$ with relations
\begin{eqnarray}
  wxw^{-1}={}^{w}x,\label{def1gr}\\
  \left[x,y\right]=\sum_{w\in W}\Omega_w(x,y)w,\label{def2gr}
\end{eqnarray}
for $x, y\in V$ and $w \in W.$ For each $w\in W$, $\Omega_w : V\times V \rightarrow \mC$ is an alternating $2$-form on $V$; we insist these forms satisfy the coherence conditions of \cite[(1.6), (1.7)]{RS}. The algebra
$\mathbf{H}_{gr}$ is a {\it graded Hecke algebra} for $W$ and
$\mathbf{H}_{gr}\cong S(V)\otimes\mC W$ as vector spaces (\cite[Lemma
1.5]{RS}).  In particular, there is a PBW-basis given by the elements of the set
$\{fw\}$, where $w\in W$, and $f$ runs through a
homogeneous basis of $\S$.
For $f$ in $S(V)$ we again write $|f|$ for the degree of $f$.
For $i \in \mZ_{\geq 0}$ let $\cB_{<i}$ be the span of all
PBW-basis elements of the form $fx$, where $f\in S(V)$, $x\in W$ such that
$f$ is homogeneous with $|f|< i$: this induces a $\mZ_{\geq0}$-filtration of
$\mathbf{H}_{gr}$. Moreover, the commutation relation \eqref{def2gr} shows
that
\begin{equation}
\label{commforusgr} \text{$|[f,g]| \leq |f| + |g| - 2$ for all
homogeneous $f, g\in S(V)$.}
\end{equation}
Recall that $s\in W$ is a {\it bireflection} if $\op{codim} V^s:=\op{rank}(\op{id_V}-s)=2$. We denote by $\cR$ the
set of all bireflections $s$ such that for any $g\in Z_W(s)$, the $W$-centraliser of $s$, the action of $g$
restricted to $V/V^s$ has determinant equal to one. The set $\cR$ plays an important role since $\Omega_g\not=0$
implies $g=e$ or $g\in\cR$ (\cite[Theorem 1.9]{RS}). Moreover, since $V$ is the (faithful) defining reflection
representation of $W$ and $\Omega_e\in((\wedge^2V)^*)^W$, we find $\Omega_e=0$. Hence relation \eqref{def2gr}
becomes
\begin{eqnarray}
\label{def2grb}
  [x,y]=\sum_{w\in \cR}\Omega_w(x,y)w.
\end{eqnarray}
Let $N\triangleleft W$ be the normal subgroup generated by $\cR$ and let ${{\bf H}_{gr}}(N)$ be the graded Hecke
algebra associated with $N$ whose structure is inherited from ${{\bf H}_{gr}}$. The following fact illustrates
once more that $\cR$ controls ${{\bf H}_{gr}}$: there is (\cite[Lemma 1.3]{Passman}) an isomorphism of algebras
\begin{eqnarray}
\label{crossedprod}
  {{\bf H}_{gr}}\cong {{\bf H}_{gr}}(N)*' W/N,
\end{eqnarray}
where ${{\bf H}_{gr}}(N)*' W/N$ is a crossed product algebra defined as follows. As a vector space it is just
${{\bf H}_{gr}}(N)\otimes \mC[W/N]$. To define the commutator relations between these two subspaces we fix for
each coset of $W/N$ one representative. Let $\{g_i\mid i\in J\}$ be the resulting complete system of coset
representatives for $W/N$ with $g_i\in [g_i]\in W/N$.

 Let $T(V)$ be the tensor algebra and $T(V)\ast W$ be the skew product algebra with the
  relations given by \eqref{def1gr}. Hence ${{\bf H}_{gr}}=(T(V)*W)/I$ where $I$ is given by
the relations \eqref{def2grb}. These relations also define an ideal, $I(N)$, of
$T(V)*N$ such that ${{\bf H}_{gr}}(N)=(T(V)*N)/I(N)$. If now $x=\sum_{n\in N}
v_nn\in T(V)\otimes\mC N$ then define
\begin{eqnarray}
  \label{commcross}
  [g_i]x=\sum_{n\in N} {}^{g_i}v_n\, g_i n g_i^{-1}\, [g_i].
\end{eqnarray}
Passing to the quotient, this defines the commutator relations between ${{\bf
    H}_{gr}}(N)$ and $\mC[W/N]$ in ${{\bf H}_{gr}}(N)*W/N$. One can show
    that, up to isomorphism, this algebra does not depend on the choice of
    representatives. However, with these choices, the isomorphism \eqref{crossedprod}
    is explicitly given as $fg\mapsto f\cdot g_ing_i^{-1}\cdot [g_i]$, where
$f\in S(V)$, $g=g_in\in W$, $n\in N$. Since ${\bf H}_{gr}(N)$ is preserved by conjugation by the subgroup $W$ of
${\bf H}_{gr},$ we note:

\begin{lemma}
\label{Gaction}
  Let $Z({{\bf H}_{gr}}(N))$ be the centre of ${{\bf H}_{gr}}(N)$ considered as
  a subalgebra of ${{\bf H}_{gr}}$ via the isomorphism~\eqref{crossedprod}. The
  $W$-action $g.h=ghg^{-1}$ for $g\in W$, $h\in{{\bf H}_{gr}}$ induces a
  $W$-action on $Z({{\bf H}_{gr}}(N))$.
\end{lemma}

\subsection{The central subalgebra}
In the special case (see \cite[Section 3]{RS}) where $W$ is a Weyl group and ${\bf H}_{gr}$ is Lusztig's graded
Hecke algebra (as introduced in \cite{Lusztigaffine}) the following result is well-known (\cite[Proposition
4.5]{Lusztigaffine}). We retain the notation $\{{\bf a}_i : 1 \leq i \leq |W| \}$ from Section~\ref{dualbases}.

\begin{proposition}
\label{centralsub}
  \begin{enumerate}
  \item The algebra $\mathbf{H}_{gr}$ has finite global dimension.
  \item The centre $Z({\bf H}_{gr})$ contains the subalgebra $\Zo_{gr}:=\S^W$.
  \item With the notation from the previous section, $\mathbf{H}_{gr}$ becomes
  a free $\Zo_{gr}$-module of finite rank with basis $$\cB_{\mathbf{H}_{gr}}
  \quad := \quad \{{\bf a}_{i}w : w\in W, 1\leq i\leq |W|\}.$$
  \end{enumerate}
\end{proposition}

The proof of this proposition will occupy the rest of this subsection. We start with some preparations. Note that
if $\Omega_w=0$ for all $w\in W$, then $\mathbf{H}_{gr}\cong S(V)* W$, the skew group algebra. Of course, the
proposition holds in this case. For any filtered algebra $B$ we denote
  by $\op{Gr}B$ its associated graded algebra. The following holds:

  \begin{lemma}
\label{lemmaA}  Let $e_N=\frac{1}{|N|}\sum_{w\in N}w$ and consider
  $\mathbf{H}^{sph}:=e_N\mathbf{H}_{gr}(N)e_N$, the spherical subalgebra of
  $\mathbf{H}_{gr}(N)$. The $\mZ_{\geq0}$-filtration on $\mathbf{H}_{gr}(N)$
  induces a filtration on $\mathbf{H}^{sph}$ and also on its centre such that
    \begin{enumerate}
    \item \label{ll1} $\op{Gr}\mathbf{H}^{sph}\cong S(V)^N$.
    \item \label{ll2} There is an isomorphism of algebras $\Psi:Z(\mathbf{H}_{gr}(N))\cong Z(\mathbf{H}^{sph})$, $z\mapsto ze_N$.
    \item \label{ll3} $\mathbf{H}^{sph}$ is commutative, in particular
    $Z(\mathbf{H}_{gr}(N))\cong \mathbf{H}^{sph}$.
    \item \label{ll4} $\op{Gr}Z(\mathbf{H}_{gr}(N))\cong S(V)^N$.
    \end{enumerate}
  \end{lemma}

  \begin{proof}
     There is an isomorphism $S(V)^N\rightarrow e_N(S(V)*N)e_N$ via $f\mapsto fe$,
    and $e_N(S(V)*N)e_N\cong e_N(\op{Gr}\mathbf{H}_{gr}(N))e_N\cong {\it
    Gr}(e_N\mathbf{H}_{gr}(N)e_N)=\op{Gr}\mathbf{H}^{sph}$. This proves
    \eqref{ll1}.
    Statements~\eqref{ll2} and \eqref{ll3} are analogous to \cite[Theorem 3.1]{EG} and \cite[Theorem 1.6]{EG} respectively;
    details can be found in \cite{Katrin}.
    Since $\Psi$ preserves the filtration and is surjective on each layer, the
    last statement follows from \eqref{ll3}.
  \end{proof}

Let $R=S(V)*N$. Recall that an associative graded algebra $(A,\diamond)$, with multiplication
  $\diamond$, is called {\it a graded
  deformation of $R$} if $A\cong R\otimes_\mC\mC[h]$ as graded vector spaces
  where $h$ is an indeterminant concentrated in degree one, $\diamond$ is
  $\mC[h]$-bilinear, and $r_1\diamond
  r_2\equiv r_1r_2 \mod hA$ for any $r_1, r_2\in R$, considered as a subspace
  of $A$. Put
  \begin{eqnarray*}
    A=A(V,N):=(T(V)[h]*N)/I_N,\quad I_N:=\langle\left[x,y\right]-\sum_{w\in
    \mathcal{R}}\Omega_w(x,y)wh^2:x,y\in V\rangle.
  \end{eqnarray*}
  Note that $I_N$ becomes homogeneous, hence $A$ is graded. It follows directly that
  $A$ is a graded deformation of $R$ and $A/(h-1)A=\mathbf{H}_{gr}(N)$.

  \begin{proof}[Proof of Proposition~\ref{centralsub}]
  The first statement is clear from \cite[Corollary 7.6.18(i)]{MR}, since $\mathbf{H}_{gr}$ is filtered such that
  $\op{Gr}(\mathbf{H}_{gr})\cong S(V)* W$ and the latter has finite global
  dimension. The last statement will follow as soon as we
  established the second.

Recall (from Lemma~\ref{Gaction}) that $W$ and hence $W/N$ act on the
  centre of $\mathbf{H}_{gr}(N)$. We get $\op{Gr}(Z(\mathbf{H}_{gr}(N))^{W/N})=(\op{Gr}Z(\mathbf{H}_{gr}(N)))^{W/N}=(S(V)^N)^{W/N}=S(V)^W$ by
  Lemma~\ref{lemmaA}, and $e_NAe_N$ is a commutative graded deformation of $S(V)^N;$ the proof of this is analogous to
  the proof of \cite[Theorem 1.6]{EG}, and is given in detail in  \cite{Katrin}. The infinitesimal commutative graded deformations are controlled by the
 second Harrison cohomology (\cite[Theorem 8]{Harrison}, \cite[Section
 4]{Gerst}). In our situation $B:=(e_NAe_N)^{W/N}$ is a (global)
 commutative graded deformation of $S(V)^W$. On the other hand, $W$ is a complex reflection
 group, hence $S(V)^W$ is a polynomial ring,  and so there are no non-trivial graded commutative deformations (\cite[Theorem 11]{Harrison}). Hence
 $B$ is a trivial deformation, and therefore $B/(h-1)B=S(V)^W$. On the
 other hand $B/(h-1)B=(e_N\mathbf{H}_{gr}(N)e_N)^{W/N}=(\mathbf{H}^{sph})^{W/N}$, hence $Z(\mathbf{H}_{gr}(N))^{W/N}=S(V)^W$
  by Lemma~\ref{lemmaA}.
The claim of the proposition follows then from \eqref{crossedprod} as follows: Let
   $f\in S(V)^W$, in particular $fg=gf\in \mathbf{H}_{gr}$ for any $g\in W$. Since  the centre of $\mathbf{H}_{gr}(N)$
   is given by $S(V)^N$ and $f\in S(V)^W\subset S(V)^N$, we get
   $fh=hf$ for any $h\in\mathbf{H}_{gr}(N)$, considered as a subspace of
    $\mathbf{H}_{gr}$. Hence, $f$ is in the centre of $\mathbf{H}_{gr}$.
  \end{proof}

\subsection{The centre} Although it is not needed for the results of this paper, we record here the fact that the
inclusion of $S(V)^W$ in the centre of $\mathbf{H}_{gr}$ is in fact an equality. In the special case where $W$ is
a Weyl group, this result is \cite[Proposition 4.5]{Lusztigaffine}.

\begin{theorem}
Retain the notation of Sections~\ref{GradedHecke} and \ref{centralsub}. Then $S(V)^W = Z(\mathbf{H}_{gr}).$
\end{theorem}
\begin{proof}
From Proposition~\ref{centralsub}(2) we know that $\Zo_{gr} := S(V)^W \subseteq Z:= Z(\mathbf{H}_{gr}).$ Let $F$ and $E$ be the
quotient fields of $\Zo_{gr}$ and $Z$ respectively, and let $Q$ be the (simple artinian) quotient ring of
$\mathbf{H}_{gr},$ so $F \subseteq E \subseteq Q.$ Since $\mathbf{H}_{gr}$ is a finitely generated module over the
commutative affine algebra $\Zo_{gr},$ $Z \cap F$ is a finitely generated $\Zo_{gr}-$module. Therefore, since $\Zo_{gr}$ is
integrally closed, $Z \cap F = \Zo_{gr}.$ Suppose for a contradiction that $\Zo_{gr} \subsetneq Z.$ Then $F \subsetneq E.$
It follows that $$ \mathrm{dim}_E(Q) < \mathrm{dim}_F(Q) = |W|^2.$$ That is, the PI-degree of $\mathbf{H}_{gr}$ is
strictly less than $|W|$, or - equivalently - the maximal dimension of an irreducible $\mathbf{H}_{gr}-$module is
strictly less than $|W|,$ \cite[Theorem I.13.5 and Lemma III.1.2]{BGo}.

We now claim that the maximal dimension of irreducible $\mathbf{H}_{gr}-$modules is $|W|.$ To see this, consider
the algebra $\hat{\mathbf{H}}_{gr}$, which has the same generators as $\mathbf{H}_{gr},$
but is constructed as an algebra over a polynomial algebra $\mathbb{C} [h].$ Relations (4.1) are unchanged, but
the right hand sides of the relations (4.2) are multiplied by $h^2$. Thus $\hat{\mathbf{H}}_{gr}$ is
$\mathbb{N}-$graded, with $h$ and the elements of $V$ having degree 1, and elements of $W$ degree 0. As before, we
can show that $\mathbb{C}[h]S(V)^W \subseteq Z(\hat{\mathbf{H}}_{gr}),$ so that $\hat{\mathbf{H}}_{gr}$ has
PI-degree at most $|W|$ by the same argument as above. On the other hand,
$\hat{\mathbf{H}}_{gr}/h\hat{\mathbf{H}}_{gr} \cong S(V)*W,$ the skew group algebra, and this has irreducible
modules of dimension $|W|$ - for example, one has an irreducible $S(V)*W$-structure on $S(V)/\mathfrak{m}S(V)$ for
any maximal ideal $\mathfrak{m}$ of $S(V)^W$ contained in a maximal orbit (of size $|W|$) of maximal ideals of
$S(V).$ Therefore,
$$ \mathrm{PI-degree}(\hat{\mathbf{H}}_{gr}) = \mathrm{PI-degree}(S(V)*W) =|W|.$$
Now the Azumaya locus of $\hat{\mathbf{H}}_{gr}$ is dense in $\mathrm{maxspec}(Z(\hat{\mathbf{H}}_{gr}))$ (\cite[Theorem III.1.7]{BGo}); in
particular, there must be an irreducible $\hat{\mathbf{H}}_{gr}-$module $U$ annihilated by $h - \lambda$ for some
$0 \not= \lambda \in \mathbb{C}.$ This implies that
$\mathrm{PI-degree}(\hat{\mathbf{H}}_{gr}/(h-\lambda)\hat{\mathbf{H}}_{gr}) = |W|,$ and so proves our claim, since
all such factors, for $\lambda \not= 0,$ are isomorphic to $\mathbf{H}_{gr}.$ We have thus obtained the desired
contradiction, so the proof is complete.
\end{proof}

\subsection{The bilinear form}
 Consider the $\Zo_{gr}$-linear map
\begin{eqnarray*}
  \Phi_{gr}:\quad\quad\quad\mathbf{H}_{gr}&\longrightarrow&\Zo_{gr}\\
  \cB_{\mathbf{H}_{gr}}\ni {\bf a}_{i}w &\longmapsto&
  \begin{cases}
1&\text{if $w=e$, $i=max$},\\
0&\text{otherwise}.
  \end{cases}
\end{eqnarray*}

Define the form $\B$ for $\mathbf{H}_{gr}$ by $\B(a,b) = \Phi_{gr}(ab),$ for $a,b \in \mathbf{H}_{gr}.$ We can now deduce the
\label{6.3gr}
\begin{lemma}
\label{maingr} The functional $\Phi_{gr}$ above satisfies Hypothesis
\ref{hyp}.
\end{lemma}
\begin{proof}
 The proof is completely analogous to Lemma~\ref{6.3}.
\end{proof}

\begin{theorem}
\label{gradedHecke}
  The graded Hecke algebra $\mathbf{H}_{gr}$ is a free Frobenius extension of its
  central subalgebra $\Zo_{gr}$ with Nakayama automorphism $\nu$
  given by $\nu(w)=\epsilon_{V}(w)^{-1}w$, $\nu(v)=v$ for $w\in W$, $v\in V$.
\end{theorem}

\begin{proof}
It is immediate from Lemmas \ref{6.3} and 2.3 that ${\bf H}_{gr}$ is a free Frobenius extension of $\Zo_{gr}$ with
form $\B$ as defined above. Therefore it remains only to determine the Nakayama automorphism. Let $\nu$ be as in
the theorem, and let $fw$ be a typical element from
  $\cB_{\mathbf{H}_{gr}}$. First, let $x\in W.$ Then
  $\B(fw,x)=\Phi_{gr}(fwx)=\delta_{w,x^{-1}}\Phi_{gr}(f)$ from the definition of
  $\Phi_{gr}$, and
  $\B(\nu(x),fw)=\Phi_{gr}(\nu(x)fw)=\Phi_{gr}(\epsilon_V(x)^{-1}xfw)=\Phi_{gr}(fxw)=\delta_{w,x^{-1}}\Phi_{gr}(f)$ using the defining relations~\eqref{def1gr} of $\mathbf{H}_{gr}$ and again the definition of $\Phi_{gr}$. If $a\in V$ we get
\begin{eqnarray}
    \B(fw,a)&=&\Phi_{gr}(fwa)=\Phi(f \, {}^wa \, w)\stackrel{(*)}{=}\Phi(faw)\stackrel{(**)}{=}\Phi(afw)=\B(a,fw).
  \end{eqnarray}
The equality (**) arises since the degree of $fa$ and
$af$ is $|f| + 1$ and so both sides are zero unless $|f|\geq N-1$. In the case $|f| = N$ or $N-1$ then $|[a,f]| < N$
by \eqref{commforusgr}. This then means that $\Phi(faw) = \Phi(
afw - [f,a]w) = \Phi (afw)$, as required. The equality (*) is true, because we have zero on both sides if
$w\not=e$. Hence $\B(fw,a)=\B(a,fw)$ holds.

Since $\mathbf{H}_{gr}$ is generated by $V$ and $W$, $\B(x,y)=\B(\nu(y),x)$ for any $x, y\in\mathbf{H}$, where
$\nu$ is as claimed.
\end{proof}

Just as in Section \ref{conseq}, we can immediately deduce the

\begin{corollary} The factor ${\mathbf{H}_{gr}}_{\chi}$ of the
graded Hecke algebra $\mathbf{H}_{gr}$ by a maximal ideal
$\mathfrak{m}_{\chi}$ of its central subalgebra $\mathcal{Z}_{gr}$
is a Frobenius algebra which in general is not symmetric.
\end{corollary}

\section{The extended affine Hecke algebra}
In this section we show that the extended affine Hecke algebra
$\mathcal{H}$ is a Frobenius
extension of its centre, with non-trivial Nakayama automorphism,
so that the corresponding reduced  algebras
$\mathcal{H}_\chi$ are Frobenius but not, in general, symmetric.
\subsection{}
Let $W$ be a (finite) Weyl group with length function $l$ and integral weight lattice $X$, and let $v$ be an
indeterminant. For a parameter set $L$ we denote by $\mathcal{H}$ the corresponding extended affine Hecke algebra
over $\mC[v,v^{-1}]$ as defined in \cite[3.1]{Lusztigaffine}. With the notation from \cite[Lemma
3.4]{Lusztigaffine} $\cH$ is a free $\mC[v,v^{-1}]$-module with basis $T_w\theta_x$, for $w\in W$, $x\in X$, and
the subalgebra $\mathbb{C}[v,v^{-1}]\langle \theta_x :x \in X \rangle$ is a Laurent polynomial algebra. Let
$\Zo_\cH=\mC[v,v^{-1}][X]^W$ be the centre of $\cH$ \cite[Proposition 3.11]{Lusztigaffine}. By the
Pittie-Steinberg Theorem (\cite{Steinberg}), $\Zo_\cH$ is a polynomial ring over $\mC[v,v^{-1}]$ and $\cH$ is free
over $\Zo_\cH$ of finite rank $|W|^2$. By abuse of language we denote by $(\mC[v,v^{-1}][X]^W_+)$ the augmentation
ideal in $\Zo_\cH$. We consider the coinvariant $\mC[v,v^{-1}]$-algebra $\mC[v,v^{-1}][X]/(\mC[v,v^{-1}][X]^W_+)$
which we equip with a $\mZ$-grading. This induces a $\mZ_{\geq 0}$-filtration on
$\mC[v,v^{-1}][X]/(\mC[v,v^{-1}][X]^W_+)$. We fix again a pair of (homogeneous) dual bases $\{\overline{{\bf
a}}_i:1\leq i\leq |W|\}$, $\{{\overline{\bf a}}^i:1\leq i\leq |W|\}$ of the coinvariant algebra and lift these
elements to bases $\{{\bf a}_i:1\leq i\leq |W|\}$, $\{{\bf
  a}^i:1\leq i\leq |W|\}$ of the free $\Zo$-module $\mC[v,v^{-1}][X]$ such
that the (filtered) degree of  ${\bf a}_i$ agrees with the grading
degree of $\overline{\bf a}_i$. Then $\cH$ is free
over $\Zo_\cH$ of rank $|W|^2$. Let $\cB_\cH$ be the basis given by the $T_w {\bf a}_i$.

\begin{lemma}
  Let $\cH_i$ be the $\Zo_\cH$-span of all $T_w{\bf a}_j$, where $1\leq
  j\leq |W|$ and $l(w)\leq i$. Then $\cH=\bigcup_{i\geq 0}\cH_i$ is a
  filtration of $\cH$.
\end{lemma}
\begin{proof}
  We have to show that $\cH_i\cH_j\subseteq\cH_{i+j}$ for any $i,
  j\in\mZ_{\geq 0}$. With the notation from \cite[Proposition
  3.9]{Lusztigaffine} we have $\theta_xT_s\equiv T_s\theta_{s(x)}\mod\cH_0$, and then for any $w\in W$
  \begin{eqnarray}
    \label{eq:Lu}
    \theta_xT_w\equiv T_w\theta_{w^{-1}(x)}\mod \cH_{l(w)-1}
  \end{eqnarray}
by induction. To establish the lemma we only have to show that
$T_w\theta_x T_v\theta_y\in\cH_{l(w)+l(v)}$ for any $v,w\in W$,
$x,y\in X$. This is of course true if $l(v)=0$. From
formula~\eqref{eq:Lu} we get $T_w\theta_x T_v\theta_y\equiv
T_wT_v\theta_x\theta_y$ modulo
$T_w\cH_{l(v)-1}\subseteq\cH_{l(w)+l(v)-1}\subset\cH_{l(w)+l(v)}$.
On the other hand $T_wT_v\theta_x\theta_y\in\cH_{l(w)+l(v)}$ and
we are done.
\end{proof}
\subsection{}
Analogous to the cases above we define a $\Zo_\cH$-linear map
\begin{eqnarray*}
  \Phi_\cH:\quad\quad\quad\cH&\longrightarrow&\Zo_\cH\\
  \cB_{\cH}\ni T_w{\bf a}_i&\longmapsto&
  \begin{cases}
1&\text{if $w=e$ and $i=max$},\\
0&\text{otherwise}.
  \end{cases}
\end{eqnarray*}

\label{6.3ext}
\begin{proposition}
\label{mainext} The functional $\Phi_\cH$ defined above satisfies
Hypothesis \ref{hyp}.
\end{proposition}

To prove this statement we need the following easily verified
formulas:

\begin{lemma} Let $w, x\in W$, $f,g\in\mC[v,v^{-1}][X]$.
  \begin{enumerate}
  \item \label{TxTw} Let $T_xT_w=\sum_{y\in W} h_y T_y$ in $\cH$. If $h_e\not=0$ then
  $w=x^{-1}$.
  \item \label{stupid} If $l(w)\geq l(x)$ then $\Phi_\cH(T_wfT_xg)\not=0$ implies $x=w^{-1}$.
  \end{enumerate}
\end{lemma}

\begin{proof}
Statement~\eqref{TxTw} is an easy induction argument using the
defining relations of $\cH$ and therefore omitted. (For a
representation theoretic interpretation of this statement we refer
to \cite[Theorem 3.1]{Stcompf}). To verify
Statement~\eqref{stupid} note that if $x\in W$, $l(w)\geq l(x)$
then there exists some $h\in\mC[v,v^{-1}][X]$ such that
  $T_wfT_x=T_wT_xh$ modulo $T_w\cH_{l(x)-1}$ (by
  formula~\eqref{eq:Lu}). Therefore we get $T_wfT_xg=T_wT_xhg+r$, where
  $r\in T_w\cH_{l(x)-1}$. Since $l(x)-1<l(w)$, using Statement~\eqref{TxTw} we
  deduce that $\Phi_\cH(r)=0$ and so
  $\Phi_\cH(T_wfT_xg)=\Phi_\cH(T_wT_xhg)$. The claim follows by applying
  Statement~\eqref{TxTw} again.
\end{proof}

\begin{proof}[Proof of the proposition]
  Let $0 \not= u \in\cH$, $u=\sum_{w,i}z_{w,i}T_w{\bf a}_i$, where
  $z_{w,i}\in\Zo_\cH$. Choose $x$ of minimal length such that
  $z_{x^{-1},i}\not=0$ for some $i$. From the lemma above and formula~\eqref{eq:Lu} we get
  \begin{eqnarray*}
  \Phi_\cH(uT_xf)=\Phi(\sum_{w,i}z_{w,i}T_w{\bf
  a}_iT_xf)=\Phi(\sum_{i}z_{x^{-1},i}T_{x^{-1}}{\bf a}_iT_xf)
  \end{eqnarray*}
  for any $f\in\mC[v,v^{-1}][X]$. Using again the lemma above and
  formula~\eqref{eq:Lu} we can rewrite
  the expression $\sum_{i}z_{x^{-1},i}T_{x^{-1}}{\bf a}_iT_x$ in the form
  $\sum_{i}c_{i} {\bf a}_i +r$, where $r\in \cH$ is such that when expanded in
  the standard bases no $T_e$ occurs, and $c_i\in\Zo_\cH$ are not
  all zero. Since $\Phi_\cH(rf)=0$ for any $f\in\mC[v,v^{-1}][X]$, it is
  enough to verify the Hypothesis~\ref{hyp} for $u=\sum_{i}c_{i} {\bf
  a}_i$. But now we are in a familiar situation, except that we have only
  filtered algebras instead of graded algebras. Nevertheless, the statement follows as in Lemma~\ref{6.3}.
\end{proof}

\begin{theorem}
\label{affineHecke}
  The extended affine Hecke algebra $\cH$ is a free Frobenius extension of its
  centre $\Zo_\cH$. In general, this extension is not symmetric.
\end{theorem}

\begin{proof}
  We only have to verify that the Nakayama automorphism is non-trivial in
  general. This however follows directly from
\cite[Theorem 9.3]{Lusztigaffine} and Theorem~\ref{gradedHecke}.
\end{proof}
\subsection{}
Just as in Section \ref{conseq}, we deduce the

\begin{corollary} The factor $\cH_{\chi}$ of the extended affine Hecke
algebra $\cH$ by a maximal ideal $\mathfrak{m}_{\chi}$ of the
centre $\mathcal{Z}_{\cH}$ is a Frobenius algebra; in general it
is not symmetric.
\end{corollary}

\subsection{Nil-Hecke algebras}
We would like to mention at least two related algebras, where our approach works, namely the  {\it affine Nil-Hecke algebra} $\cH^{nil}$ and the {\it graded affine Nil-Hecke algebra} $\cH^{nil}_{gr}$ associated to a Weyl group $W$. (For the definitions see e.g. \cite{GR}). Analogous to the affine Hecke algebra case, the centre of $\cH^{nil}$ is $\Zo=\mC[X]^W$ and $\cH^{nil}$ is a free $\Zo$-module of rank $|W|^2$ (\cite[(1.9)]{GR}), similarly for the  graded affine Nil-Hecke algebras. If we define the forms completely analogous to the affine and graded Hecke algebras we deduce that $\cH^{nil}$ and $\cH^{nil}_{gr}$ are free Frobenius extensions over their centres.

\section{The quantised universal enveloping algebra}
In this section we show that the quantised enveloping algebra
$\cU_\epsilon(\mg)$ at a root of unity $\epsilon$ is a Frobenius
extension of its Hopf centre, with trivial Nakayama automorphism,
so that the reduced quantised enveloping algebras
$\cU_\epsilon(\mg)_\chi$ are symmetric.
\subsection{The PBW-basis and the central subalgebra}
Let $\mg$ be a complex semisimple Lie algebra. We fix a Borel and
Cartan subalgebra  of $\mg$, $\mb\supseteq\mh ,$ and denote the
Weyl group by $W$ and the set of simple reflections by $S$. Let
$\pi$ be the corresponding set of simple roots and $\rho$ the
half-sum of positive roots. Let $\epsilon\in\mC$ be an $l$-th root
of unity, for some odd positive integer $l$, $l\not=3$ if $\mg$
has a summand of type $G_2$. Let $Q\subseteq P$ be, respectively,
the root lattice and the  weight lattice of $\mg$, with the
$W$-equivariant bilinear form $(\,,\,):P\times Q\rightarrow \mZ$.

The simply connected form of the quantised universal enveloping algebra $\cU=\cU_\epsilon(\mg)$ is a $\mC$-algebra
with generators $E_\alpha$, $F_\alpha$, $K_\la$, for $\alpha\in\pi$ and $\la\in P$. For the defining relations and
further details we refer for example to \cite[9.1]{ChariPresley} or \cite[I.6.3, III.6.1]{BGo}. Let $w_0$ be the
longest element of $W,$ and fix a reduced expression
\begin{eqnarray}
\label{w0}
  w_0=s_{i_1}s_{i_2}\ldots s_{i_N},
\end{eqnarray}
where $s_{i_j}\in S$ for $1\leq j\leq N$. Let $\alpha_{i_j}$ be
the simple root corresponding to  $s_{i_j}\in S$. Recall that
Lusztig defined an action on $\cU$ of the braid group $B$
corresponding to $W,$ (see \cite{Luquantroot}, \cite[Section
9]{ChariPresley}, \cite[Section 8]{Jquant} or \cite[I.6.7, I.6.8]{BGo}). Let
$T_i$ be the automorphism in $B$ corresponding to the simple
reflection $s_i\in S$. We set
\begin{eqnarray}
  \label{eq:beta}
\beta_k:=s_{i_1}s_{i_2}\ldots s_{i_{k-1}}(\alpha_{i_k}),
\end{eqnarray}
and put $E_{\beta_k}=T_{i_1}T_{i_2}\ldots
T_{i_{k-1}}(E_{\alpha_{i_k}})$ and
$F_{\beta_k}=T_{i_1}T_{i_2}\ldots T_{i_{k-1}}(F_{\alpha_{i_k}})$.
For any sequence $\mathbf{m}=(m_1,m_2,\ldots,
m_N)\in\mathbb{Z}^N_{\geq0}$ let
  \begin{eqnarray*}
E^\mathbf{m}&=&E_{\beta_1}^{m_1}E_{\beta_2}^{m_2}\ldots E_{\beta_N}^{m_N},\\
F^\mathbf{m}&=&F_{\beta_N}^{m_N}F_{\beta_{N-1}}^{m_{N-1}}\ldots
  F_{\beta_1}^{m_1}.
  \end{eqnarray*}
   This yields a PBW-basis of $\cU$ (associated with \eqref{w0}),
   namely
  \begin{eqnarray*}
   \cB \quad = \quad \{  F^{\bf k} K_\la E^{\bf m} : {\bf k}, {\bf m}\in\mZ^{N}_{\geq0},\la\in P \},
  \end{eqnarray*}
see \cite[Theorem 9.3]{ChariPresley}, \cite[I.6.2, III.6.1]{BGo}.
The subspace $\Zo$ of $\cU$ spanned by the monomials $F^{l\bf k}
K_{l \la} E^{l \bf m}$ is a central Hopf subalgebra  of $\cU$,
called the \emph{$l$-centre}, and $\cU$ is a free $\Zo$-module of
finite rank (see \cite[19.1]{ChariPresley}, \cite[III.6.2]{BGo}).
As a $\Zo$-basis of $\cU$ one can choose the subset $\cB'$ of
$\cB$ given by elements of the form
\begin{eqnarray}
\label{B'} F^{\bf k} K_\la E^{\bf m},
  \end{eqnarray}
where $0\leq k_i, l_i< l$ and the coefficients of $\la$ in terms
of fundamental weights are non-negative integers less than $l$.

\subsection{Filtrations, degrees and commutation
formulas}\label{filt} To simplify formulas we set
$E_{i}=E_{\beta_i}$ and $F_{i}=F_{\beta_i}$. (Note that $E_i$ is
not $E_{\alpha_i}$ in general.) Let $i<j$. There are commutation
formulas holding in $\cU$ as follows \cite[Proposition I.6.10, Theorem III.6.1(4)]{BGo}:
\begin{eqnarray}
\label{Es}
  E_{i}E_{j}&=&\epsilon^{(\beta_i,\beta_j)}E_{j}E_{i}+r\\
\label{Fs}
  F_{i}F_{j}&=&\epsilon^{-(\beta_i,\beta_j)}F_{j}F_{i}+r'
\end{eqnarray}
where $r$ (resp. $r'$), written in the PBW-basis, involves no
monomial containing any $E_{k}$ (resp. $F_{k}$) for $k\leq i$ or
$k\geq j$.

The algebra $\cU$ is $Q$-graded (see e.g. \cite[4.7]{Jquant}), but
 also has
  several other filtrations, \cite[10.1]{ChariPresley}, \cite[I.6.11, III.6.1]{BGo}. First, there is
  the \emph{degree filtration}, a $\mZ_{\geq
  0}$-filtration obtained by putting $F^{\bf k} K_\la E^{\bf m}\in\cB$ in degree
\begin{eqnarray*}
  \op{deg}(F^{\bf k} K_\la E^{\bf m})=\sum_{i=1}^N(k_i+m_i)\op{ht}(\beta_i),
\end{eqnarray*}
where $\op{ht}$ denotes the height function. One can refine this
to a $(\mZ_{\geq 0})^{2N+1}$-filtration by putting $F^{\bf k}
K_\la E^{\bf m}\in\cB$ in degree
\begin{eqnarray*}
  \op{d}\big(F^{\bf k} K_\la E^{\bf m}\big)=\big(k_N, k_{N-1},\ldots k_1,m_1,m_2\ldots
  m_N, \op{deg}(F^{\bf k} K_\la E^{\bf m})\big).
\end{eqnarray*}
Putting the reverse lexicographic ordering on $(\mZ_{\geq
  0})^{2N+1}$ (\ie $e_1<e_2<\ldots$, where $(e_i)_j=\delta_{i,j}$) defines the filtration by \emph{total degree}.
 The $E$'s and $F$'s commute up to terms of
lower total degree, \cite[10.1]{ChariPresley}, \cite[Proposition I.6.11]{BGo}:
  \begin{eqnarray}
  \label{EF}
    E_{i}F_j=F_jE_i+\text{  terms of lower total degree}.
  \end{eqnarray}
We denote by $$\op{max}\quad := \quad 2(l-1)\sum_{i=1}^N
\op{ht}(\beta_i)$$ the maximal
  $\op{deg}$-value on $\cB'$.

\subsection{The bilinear form}\label{form}
  In view of the $\Zo$-freeness of $\cU$ on the basis $\cB',$ we can define a $\Zo$-linear map $\Phi:\cU\rightarrow
  Z$ as follows. Set $\mathbf{l} := (l-1,l-1,\ldots l-1),$ and
  define
  \begin{eqnarray*}
      \Phi: \cB' \longrightarrow \Zo : F^{\bf k} K_\la E^{\bf m}&\longmapsto&
      \begin{cases}
        1&\text{ if $\mathbf{k}=\mathbf{m}=:\mathbf{l}$, $\la=0$},\\
        0&\text{ otherwise,}
      \end{cases}
  \end{eqnarray*}
and extend this $\Zo$-linearly.
  \begin{lemma}
\label{nondeg}
The functional $\Phi$ satisfies Hypothesis~\ref{hyp}.
  \end{lemma}
  \begin{proof}
    For $\mathbf{m}\in(\mZ_{\geq0})^N$, define $\mathbf{\tilde m}:=\mathbf{l}-\mathbf{m}\in(\mZ)^N$. For $x= F^{\bf
    k} K_\la E^{\bf m}\in \cB$ and $\mu \in P$ set $\tilde{x}_\mu= F^{\bf \tilde k} K_\mu
    E^{\bf\tilde m}$ and write $k_i(x)=k_i$, $m_i(x)=m_i$.\\

 {\it Claim 1: Let $x$, $y\in\cB$. If $\op{deg}(x)+\op{deg}(y)<\op{max}$ then
    $\Phi(xy)=0$.} \\
This follows directly from the fact that the commutation relations
(\ref{Es}), (\ref{Fs}) and (\ref{EF}) do not increase
    the $\op{deg}$-value and $\Phi$ annihilates every monomial in $\cB'$ which is not of maximal
    $\op{deg}$-value.\\

 {\it Claim 2: Let $x, y\in\cB'$, $\mu\in P$. If $\op{d}(x)<\op{d}(\tilde{y}_\mu)$ then
    $\Phi(yx)=0$.}  \\
 If $\op{d}(x)<\op{d}(\tilde{y}_\mu)$ then
    $\op{deg}(x)\leq\op{deg}(\tilde{y}_\mu)=\op{max}-\op{deg}(y)$, hence
    $\op{deg}(x)+\op{deg}(y)\leq\op{max}$. By Claim 1 we only have to deal
    with the case $\op{deg}(x)+\op{deg}(y)=\op{max}$. From our assumption and
    the definition of $\op{d}$ it follows that \emph{either}
    \begin{itemize}
    \item there is a $k_j(x)$ such that
     $k_j(x)\not=l-1-k_j(y)$ (so that $k_j(x)<l-1-k_j(y)$), \newline
\noindent \emph{or}
    \item there is  an
    $m_j(x)$ such that $m_j(x)<l-1-m_j(y)$ (so that $m_j(x)<l-1-m_j(y)$).
    \end{itemize}
    Let us assume the latter for a
    moment and choose $j$ maximal with this property. Recalling from
    \eqref{EF} and the commutation relations for the $K$s with the $E$s and the $K$s with the $F$s
    that the relevant generators of $\cU$
    commute up to nonzero scalars and terms of lower
    $\op{deg}$-value, we see that it is enough to show that
    \begin{eqnarray*}
      \Phi(F^{{\bf k}(y)}F^{{\bf k}(x)}K_{\la(x)}K_{\la(y)}E^{{\bf m}(y)}E^{{\bf m}(x)})=0.
    \end{eqnarray*}From the relation~\eqref{Es} we get
    \begin{eqnarray*}
    E_N^{m_n(y)}E_1^{m_1(x)}\cdots E_N^{m_N(x)}&=&cE_1^{m_1(x)}\cdots E_N^{m_N(x)+m_N(y)}+r,
    \end{eqnarray*}
for some $c\in\mC^*$ and some $r \in \cU$ such that $E_N$ occurs
in every monomial in $r$ with power strictly smaller than
$m_N(y)+m_N(x)\leq l-1$. In particular,
\begin{eqnarray*}
 && \Phi(F^{{\bf k}(y)}F^{{\bf k}(x)}K_{\la(x)}K_{\la(y)}E^{{\bf m}(y)}E^{{\bf
  }{\bf{m}}(x)})\\
&=&c\Phi(F^{{\bf k}(y)}F^{{\bf k}(x)}K_{\la(x)}K_{\la(y)}E_1^{m_1(y)}\cdots E_{N-1}^{m_{N-1}(y)}E_1^{m_1(x)}\cdots
  E_N^{m_N(x)+m_N(y)}).
\end{eqnarray*}
Repeating this argument we get
\begin{eqnarray*}
&&  \Phi(F^{{\bf k}(y)}F^{{\bf k}(x)}K_{\la(x)}K_{\la(y)}E^{{\bf m}(y)}E^{{\bf m
  }(x)})=c'\Phi(F^{{\bf k}(y)}F^{{\bf k}(x)}K_{\la(x)+\la(y)}X),
\end{eqnarray*}
where $X=E_1^{m_1(y)}\cdots
  E_{j-1}^{m_{j-1}(y)}E_1^{m_1(x)}\cdots
  E_{j-1}^{m_{j-1}(x)}E_j^{m_j(x)+m_j(y)}E_N^{m_N(x)+m_N(y)}$ and $c' \in \mC^*$.
The result is zero since any commutation of the $E_i$s for $i<j$
does not involve $E_j$ because of \eqref{Es}, and since
$m_j(x)+m_j(y)<l-1$. The remaining case, where there is a $k_j(x)$
such that
     $k_j(x)\not=l-1-k_j(y),$ can be proved similarly and is therefore omitted. \\

{\it Claim 3: Let $x$, $y\in\cB$, $x=F^{{\mathbf k}}K_\la E^{\mathbf{m}}$,
  $\mu \in P$ and $\op{d}(x)=\op{d}(\tilde{y}_\mu)$ then we have $\Phi(yx)\not=0$ if and only if $\la+\mu\in lP$.}\\
With the arguments from the proof of Claim 2 we get
\begin{eqnarray*}
  \Phi(yx)=c\Phi(F^{\mathbf{l}} K_{\la(x)+\la(y)} E^\mathbf{l}),
\end{eqnarray*}
for some nonzero number $c\in\mC$. Claim 3 follows then from the definition
of $\Phi$.\\

To prove the proposition, let $x\in\cU$ be arbitrary and write
$x=\sum_{y\in\cB'} z_y y$ with $z_y\in \Zo$. We choose
$b=F^{{\mathbf k}}K_\la E^{\mathbf{m}}\in\cB'$ of maximal total
degree such that $z_{b}\not=0$. If now $a=F^{{\mathbf r}}K_\nu
E^{\mathbf{s}}\in\cB'$ with $z_a\not=0$ and $\nu\in P$ arbitrary,
then we have have $\op{deg}(a)\leq\op{deg}(b)$, hence,  for any
$\mu\in P,$
$$\op{deg}(a)+\op{deg}(\tilde{b}_\mu)=\op{deg}(a)+\op{max}-\op{deg}(b)\leq\op{max}.$$
 If this inequality is strict, Claim 1 implies
$\Phi(\tilde{b}_\mu a)=0$. Let us assume equality. Then we either
have $\op{d}(a)<\op{d}(b)$ which implies $\Phi(\tilde{b}_\mu a)=0$
by Claim 2, or $\op{d}(a)=\op{d}(b)$. The latter means (because of
Claim 3) that $\Phi(\tilde{b}_{\mu} a)=0$ except when $\nu+\mu\in
lP$. In particular, $\Phi(\tilde{b}_{-\la}a)=0$, except when $a=b$.

Summarising, we get
$\Phi(\tilde{b}_{-\la}x)=\Phi(\tilde{b}_{-\la}b)=cz_b\not=0$ for
some unit $c$, as required.
\end{proof}

\subsection{Symmetry of the form}
\label{symm}
\begin{lemma}
\label{NakU}
  The Nakayama automorphism $\nu$ of $\cU$ with respect to $\mathbb{B}$ is the identity.
\end{lemma}

\begin{proof}
  We have to prove that $\nu$ fixes all generators. We will run through all possibilities $y$ for generators and prove $\B(x,y) = \B(y,x)$ for all $x\in \mathcal{B}'$.

  First, let $y=K_\la$. From
  Claim 2 we have automatically $\B(x,y)=0=\B(y,x)$ unless $x$ has maximal
  degree $\op{max}$. But then $K_\la$ commutes with $x$ and hence
  $\nu(K_\la)=K_\la$.

Now let $y=E_{\alpha}$ for some simple root $\alpha$, and let $x \in \cB'.$ Claim~1 implies that
  $\B(x,y)=0=\B(y,x)$ unless $\deg{x}\geq\op{max}-1$, because
  $\op{deg}(E_{\alpha})=1$. If $x = F^{\bf l}K_{\lambda} E^{\bf l}$ then both $yx$ and $xy$ have $Q$-grade equal to $\alpha$. Thus $\Phi(yx) = 0 = \Phi(xy)$ since, by definition, $\Phi$ is non-zero only on elements whose $Q$-grade belongs to $\ell Q$. We now have two possibilities for $x$: either
  $\op{deg}(F^{{\bf k}(x)})\not=\frac{\op{max}}{2}$ and
  $\op{deg}(E^{{\bf m}(x)})=\frac{\op{max}}{2},$ or vice versa. Let us consider
  the first case. Then $\B(x,y)=0=\B(y,x)$, because $\Phi$
  annihilates everything which does not have the same $Q/lQ$-grading as
  $F^{\mathbf{l}}E^{\mathbf{l}}$ by definition. In the second case, the
  $Q$-grading again implies $\B(x,y)=0=\B(y,x)$, unless $m_j(x)\not=l-1$ implies
  $\beta_j=\alpha$. Let $j$ be such that this equation holds. That means we have to compare $\B(x,y)=\Phi(xE_{\alpha})$ and
  $\B(y,x)=\Phi(E_{\alpha}x)$, where $x=F^{\mathbf{l}}K_\la E_1^{l-1}\ldots
  E_{j-1}^{l-1}E_{j}^{l-2}E_{j+1}^{l-1}\ldots E_N^{l-1}$. Both terms are
  trivial unless $\la=0$.  From the commutator relation~\eqref{Es} it follows
  that $\B(x,y)=\epsilon^{-(l-1)(\beta_{j+1}+\ldots+\beta_N, \beta_j)}$ and
  $\B(y,x)=\epsilon^{-(l-1)(\beta_{1}+\ldots+\beta_{j-1}, \beta_j)}$. It is
  now enough to show that the exponents are the same.

Put $w=s_{i_1}s_{i_2}\cdots s_{i_{j-1}}$. Then $M^-=\{\beta_r :
1\leq r\leq j-1\}$ (resp. $M^+=\{\beta_r : j\leq r\leq N\}$) is
exactly the set of all positive roots such that $w^{-1}(\beta)$ is
negative (resp. positive). Set $M_1=w^{-1}(M^+)$ and
$M_2=-w^{-1}(M^-)$. The disjoint union of these two sets is
exactly the set of all positive roots (see e.g. \cite[I.4.3,
Theorem B]{Kane}). By definition (see \eqref{eq:beta}) we have
$w^{-1}(\beta_{j})=\alpha_{i_j}\in M_1$. From the definition of
$\rho$, the half-sum of positive roots, we get
\begin{eqnarray*}
(\alpha_{i_j},\alpha_{i_j})&=&(\alpha_{i_j},\alpha_{i_j}) +\sum_{\beta\in
    M_1\backslash\{\alpha_{i_j}\}}
(\beta,\alpha_{i_j})+\sum_{\beta\in M_2}(\beta,\alpha_{i_j}).
\end{eqnarray*}
Since the bracket $(\,,\,)$ is non-degenerate and $W$-equivariant, we get
\begin{eqnarray*}
  0&=&\sum_{\beta\in M_1\backslash\{\alpha_{i_j}\}}(w(\beta),w(\alpha_{i_j}))+\sum_{\beta\in
    M_2}(w(\beta),w(\alpha_{i_j}))\\&=&\sum_{\beta\in
    M^+\backslash\{\beta_{j}\}}(\beta,\beta_{j})-\sum_{\beta\in
    M^-}(\beta,\beta_{j}).
\end{eqnarray*}
Hence we get the required equality for the exponents and therefore $\B(x,y)=\B(y,x)$. We are left with the case
$y=F_\alpha$ for some simple root $\alpha$. The arguments there are similar, and therefore omitted. This completes
the proof of the lemma.
\end{proof}

\subsection{}
\label{symmetricU}

Recall the terminology of ~\ref{frob}. From Proposition
\ref{nondeg} together with Proposition \ref{cruc} and Lemma
\ref{NakU}, we have:

\begin{theorem}
  The quantised universal enveloping algebra $\cU=\cU_\epsilon(\mg)$ at an
  $l$-th root of unity is a free Frobenius extension of its $l$-th centre
  $\Zo$. The form $\B$ has a trivial Nakayama automorphism.
\end{theorem}

The following corollary is immediate from the theorem and the discussion in 2.4, noting that $\cU_{\epsilon}(\mg)$
has finite global dimension by \cite[Theorem 2.3]{BGoo}.
\begin{cor}
Let $\chi$ be a maximal ideal of $\Zo$. \begin{enumerate}[1.)]
\item The reduced quantised
  enveloping algebra $\cU_{\chi}:=\cU_\epsilon(\mg)/\cU_\epsilon(\mg)\chi$ is a symmetric algebra.
  \item  $\cU_\epsilon(\mg)$ is a Calabi-Yau $\Zo-$algebra of dimension $\dim \mg .$
  \end{enumerate}
\end{cor}

\section{Quantum Borels}
In this section we show that the quantum Borel $\cU^{\geq 0}$ at a root of unity $\epsilon$ is a Frobenius
extension of its Hopf centre, with non-trivial Nakayama automorphism, so that the reduced quantum Borels
$\cU^{\geq 0}_\chi$ are Frobenius, but not in general symmetric.
\subsection{}
\label{Borels}
Let $\mg$ be as above. Let $\cU_{\epsilon}^{\geq 0}$ be the subalgebra of $\cU_{\epsilon}(\mathfrak{g})$ generated by all the
$E$s and $K$s. The PBW-basis of $\cU_{\epsilon}(\mathfrak{g})$ gives rise to a PBW-basis of $\cU_{\epsilon}^{\geq 0}$ given by the elements of the
form $K_\la E^{\bf
  m}$, where ${\bf m}\in\mZ^{N}_{\geq0}$ and $\la\in P,$ \cite[9.3]{ChariPresley}. Moreover, $\cU_{\epsilon}^{\geq 0}$ is free over
  $Z_+:=\Zo\cap \cU_{\epsilon}^{\geq 0}$ with basis $\cB'_+$ given by all elements of the form $K_\la E^{\bf m}$, where $0\leq
m_i< l$ and the coefficients of $\la$ in terms of fundamental weights are non-negative integers less than $l$,
\cite[19.1]{ChariPresley}.

\subsection{The bilinear form and its Nakayama automorphism}
\label{4.1} Analogously to Section 4, we define a $Z_+$-linear map
 $\Phi_+:\cU_{\epsilon}^{\geq 0}\rightarrow Z_+$ by
  \begin{eqnarray*}
      \cB_+'\ni K_\la E^{\bf m}&\longmapsto&
      \begin{cases}
        1&\text{ if $\mathbf{m}=\mathbf{l}$, $\la=0$}\\
        0&\text{ otherwise.}
      \end{cases}
  \end{eqnarray*}
Define a $Z_+$-bilinear associative form $\B_+$ on $\cU_{\epsilon}^{\geq 0}$
by putting $\B_+(x,y)=\Phi(xy)$ for any $x$, $y\in\cU^{\geq 0}$.

\begin{theorem}  \label{NakBorels} Let $\cU_{\epsilon}^{\geq 0}$ be
the quantum Borel defined in (\ref{Borels}), with central
subalgebra $Z_+$ as defined there. Let $\Phi_+$ and $\B_+$ be as above.
\begin{enumerate}[1.)]
\item $\Phi_+$ satisfies Hypothesis 2.3. \item The form $\B_+$ is
non-degenerate and has a dual free pair of bases, so that
$\cU_{\epsilon}^{\geq 0}$ is a free
  Frobenius extension of $Z_+ .$
  \item The corresponding Na\-ka\-ya\-ma automorphism $\nu_+$ of $\cU_{\epsilon}^{\geq 0}$ is given by
  $\nu_+(E_\alpha)=E_\alpha$ for simple roots $\alpha$ and
  $\nu_+(K_\la)=\epsilon^{(2\rho,\la)}K_\la$ for $\la\in P$.
  \end{enumerate}
\end{theorem}

\begin{proof} The proofs of 1.) and 2.) are similar to, but easier
than the corresponding arguments for $\cU_\epsilon(\mg)$, so we
leave the details to the reader.

Consider now part 3.). As in the proof of Lemma \ref{NakU}, $\nu_+(E_\alpha)=E_\alpha$ for any simple root
$\alpha$. By the degree argument from the same proof, the value of $\nu_+(K_\la)$ is determined by
$E^\mathbf{l}K_\la=\nu_+(K_\la)E^{\mathbf{l}}$. Hence
$$\nu_+(K_\la)=\epsilon^{-(l-1)(\beta_1+\cdots\beta_N,\la)}K_\la=\epsilon^{(2\rho,\la)}K_\la.$$ The result follows.
\end{proof}

\begin{remarks}
\label{winding}
  {\rm 1. With the standard comultiplication of \cite[I.6]{BGo}, \cite[Chapter 4]{Jquant}, $E_{\alpha}\mapsto E_\alpha\otimes
    1+K_\alpha\otimes E_\alpha$, $F_{\alpha}\mapsto F_\alpha\otimes
    K_\alpha^{-1}+1\otimes F_\alpha$, $K_\la\mapsto K_\la\otimes K_\la$ for
    $\alpha\in\pi$, $\la\in P$, then
    $\nu_+$ is nothing else than the right winding automorphism \cite[I.9.25]{BGo}
    $\tau^r_{2\rho}$ of $\cU_{\epsilon}(\mathfrak{g})$
    associated with the representation $2\rho,$ restricted to $\cU_{\epsilon}^{\geq 0}$.

    2. Calculations parallel to the above will of course handle $\cU_{\epsilon}^{\leq 0},$ the
    Hopf subalgebra of $\cU_\epsilon(\mg)$ generated by the $F_{\alpha}$s and the
$K_{\lambda}$s. A more elegant approach is to make use of the Chevalley involution $\omega$ \cite[Lemma
4.6(a)]{Jquant}}: $\omega(E_{\alpha}) = F_{\alpha}$ and $\omega(K_i) = K_{i}^{-1},$ so $\omega$ is an algebra
automorphism and a coalgebra anti-automorphism. Thus one calculates that the Nakayama automorphism of
$\cU_{\epsilon}^{\leq 0},$ namely $\omega \circ \tau^r_{2\rho} \circ \omega^{-1},$ is the restriction of the
automorphism $\tau^{\ell}_{-2\rho}$ of $\cU_{\epsilon}(\mathfrak{g}).$

\end{remarks}

\section{Quantised function algebras}\label{qfun}
In this section we show that the quantised function algebra
$\cO_{\epsilon}[G]$ at a root of unity $\epsilon$ is a Frobenius
extension of its Hopf centre, with non-trivial Nakayama automorphism,
so that the reduced quantised function algebras
$\cO_{\epsilon}[G](g)$ are Frobenius but not, in general, symmetric.
\subsection{Preliminaries}\label{prefun} Let $G$ be the simply connected, semisimple algebraic group over
$\mC$ associated with the semisimple Lie algebra $\mg$. Let $B$ be the Borel subalgebra of $G$ associated with
$\pi$ and let $B^-$ be the opposite Borel. Let $T$ be the corresponding maximal torus. let $\epsilon$ be as in
(6.1), and let $\cO_\epsilon[G]$ be the quantised function algebra of $G$ at the root of unity $\epsilon.$ For the
definition and basic properties of $\cO_\epsilon [G]$, see \cite{CL} or \cite[III.7.1]{BGo}.\ff{Note, however,
that the algebra in \cite{CL} is the \emph{opposite algebra} to that in \cite{BGo}; put in another way, there is a
switch between $\epsilon$ and $\epsilon^{-1}$ in going from \cite[III.7.1]{BGo} to \cite{CL}.} Recall that de
Concini and Lyubashenko show \cite{CL} that $\cO_\epsilon [G]$ is a noetherian Hopf $\mC$-algebra which is a
finitely generated module over its centre. (An outline proof is also provided in \cite[Theorem III.7.2,
III.7.3]{BGo}.) Indeed, more specifically, $\cO_\epsilon [G]$ contains a copy of the coordinate ring of G, $\cO
[G]$, as a central Hopf subalgebra, and, by \cite[Proposition 2.2]{IK}, $\cO_\epsilon [G]$ is a free $\cO
[G]$-module of rank $l^{\mathrm{dim}G}$.

Calculations with $\cO_\epsilon[G]$ are most easily carried out by embedding it as a subalgebra of
$\cU_{\epsilon}^{\leq 0} \otimes \cU_{\epsilon}^{\geq 0}$, as in \cite[Section 4.3]{CL}. But in fact \cite{CL}
works with $(\cO_\epsilon [G])^{op}$, in terms of the definition of the function algebra of \cite{BGo} or
\cite{Jquant}; the simplest way to accommodate this here is to include a map from $\epsilon$ to $\epsilon^{-1}$
into the embedding. Once this is done, the inclusion $\mu''$ of \cite[4.3]{CL} is given by the composite
    \begin{eqnarray*}
      i': \cO_\epsilon[G]\stackrel{\text{comult}}\longrightarrow
      \cO_\epsilon[G]\otimes \cO_\epsilon[G]\rightarrow
      \cO_\epsilon[B]\otimes \cO_\epsilon[B^-]\longrightarrow \cU_{\epsilon^{-1}}^{\leq 0} \otimes_\mC
      \cU_{\epsilon^{-1}}^{\geq 0},
    \end{eqnarray*}
where the second map is the canonical one (given by ``restriction'') and the last map combines the isomorphism
from \cite[Lemma 3.4]{CL} with the parameter switch explained above. Note in passing that this embedding shows
that $\cO_\epsilon [G]$ is a domain. Moreover, by \cite[Theorem 4.6, Lemma 4.3 and Proposition 6.5]{CL}, there is
a nonzero element $z$ of $\cO [G]$, such that $i'$ extends to an inclusion
\begin{eqnarray*}
  i: \cO_\epsilon[G][z^{-1}]\longrightarrow \cU_{\epsilon^{-1}}^{\leq 0}\otimes_\mC \cU_{\epsilon^{-1}}^{\geq 0},
\end{eqnarray*}
with image generated by the elements $1\otimes E_\alpha$, $F_\alpha\otimes 1$ and $K_{-\la}\otimes K_{\la}$, for
simple roots $\alpha$ and integral weights $\la$. In the following we will often identify
$\cO_\epsilon[G][z^{-1}]$ with its image under $i$.

In particular, making this identification, a basis $\cB_\cO$ of $\cO_\epsilon [G][z^{-1}]$ as a free $\cO
[G][z^{-1}]$-module is given by the set of elements
\begin{eqnarray*}
  F^{\bf k} K_{-\la}\otimes K_\la E^{\bf m},
\end{eqnarray*}
where $0\leq k_i, m_i< l$ and the coefficients of $\la$ in terms of fundamental weights are non-negative integers
less than $l;$ for this, see the proof of \cite[Proposition 7.2]{CL}.

\subsection{The bilinear form}
\label{5.2}
We can define a $\cO[G][z^{-1}]$-linear map
$$\Phi:\quad\cO_\epsilon[G][z^{-1}]\times \cO_\epsilon[G][z^{-1}]\rightarrow \cO[G][z^{-1}]$$ by mapping
  \begin{eqnarray*}
      \cB_\cO\ni F^{\bf k} K_{-\la} \otimes K_\la E^{\bf m}&\longmapsto&
      \begin{cases}
        1&\text{ if $\mathbf{k}=\mathbf{m}=\mathbf{l}$, $\la=0$},\\
        0&\text{ otherwise,}
      \end{cases}
  \end{eqnarray*}
and extending $\cO[G][z^{-1}]$-linearly. Since $\cO[G][z^{-1}]$ is central, we get an associative
$\cO[G][z^{-1}]$-bilinear form $\B:\cO_\epsilon[G][z^{-1}]\times\cO_\epsilon[G][z^{-1}]\rightarrow
\cO_\epsilon[G][z^{-1}]$ by putting
$\B(x,y)=\Phi(xy)$ for $x$, $y\in\cO_\epsilon[G][z^{-1}]$.

\subsection{Frobenius extension}
\label{Frobenius} We can now record the key
\begin{lemma} The functional $\Phi$ satisfies Hypothesis 2.3.\end{lemma}
\begin{proof}  The argument is similar to the ones used to prove Lemma 6.3 and Theorem 7.2, and is
therefore left to the reader. \end{proof} As usual, the above lemma yields at once the first part of the following

\begin{theorem}\begin{enumerate}
\item $\cO_\epsilon[G][z^{-1}]$ is a free Frobenius extension of $\cO[G][z^{-1}]$ with the form $\B$ defined in
Section~\ref{5.2}.
\item In the notation of Remarks 7.2, the Nakayama automorphism of $\cO_\epsilon[G][z^{-1}]$ is the restriction of
the automorphism $\tau^{\ell}_{-2\rho} \otimes \tau^r_{2\rho}$ of $\cU_{\epsilon^{-1}}^{\leq 0} \otimes_\mC
      \cU_{\epsilon^{-1}}^{\geq 0}.$ In particular, it fixes $F_{\alpha} \otimes 1$ and $1 \otimes E_{\alpha}$ for
      all simple roots $\alpha$, and maps $K_{\lambda} \otimes K_{-\lambda}$ to $\epsilon^{2(2\rho,\lambda)}K_{\lambda} \otimes K_{-\lambda}$
\item There is a non-degenerate $\cO[G]$-bilinear form  $\B'$ on $\cO_\epsilon[G]$ with
  values in $\cO[G]$ and Nakayama automorphism $\nu_{\mathcal{O}} = \tau_{-2\rho}^l \otimes \tau_{2\rho}^r$.
\end{enumerate}
\end{theorem}

\begin{proof} (2) This is clear from Theorem 7.2 and Remarks 7.2(2).

(3) Choose a finite generating set $\mathcal{F}$ of $\cO_\epsilon[G]$ as a $\cO[G]$-module. There is a
non-negative integer $k$ such that $\B(u,v) \in z^{-k}\cO[G]$ for all $u,v \in \mathcal{F}.$ Let $k_0$ be the
minimal such integer, and define $\B' := z^{k_0} \B.$ Then $\B'$ has the stated properties.
\end{proof}

\begin{remark}
Suppose that $G = SL(n,\mathbb{C}),$ so that $\mathcal{O}_{\epsilon}[G]$ is generated by $\{X_{ij} : 1 \leq i,j
\leq n \}$, with the relations given at \cite[I.2.2,I.2.4]{BGo}. Then it is easy to calculate that the
automorphism $\nu_{\mathcal{O}}$ of the theorem is given by $\nu_{\mathcal{O}}(X_{ij}) \epsilon^{2(n+1-i-j)}X_{ij},$ for $i,j = 1, \ldots , n.$
\end{remark}

\subsection{Finite dimensional factors} \label{factors}
Corollary \ref{Frobenius} is sufficient to yield the desired applications to the finite dimensional representation
theory of $\cO_\epsilon[G],$ as follows:

\begin{theorem}
  Let $g \in G$ and let $\mathfrak{m}_g$ be the corresponding maximal ideal of $\cO[G]$. Then the
  algebra
  $\cO_\epsilon[G](g):=\cO_\epsilon[G]/\cO_\epsilon[G]\mathfrak{m}_g$ is a Frobenius
  algebra with Nakayama automorphism induced from $\nu_{\mathcal{O}}.$
\end{theorem}

\begin{proof}
  First let $\mathfrak{m}$ be a maximal ideal of the algebra  $\cO[G][z^{-1}]$ of Proposition \ref{Frobenius}.
  Then Proposition~\ref{Frobenius} implies that there is a
  non-degenerate $\mC$-bilinear form $\overline{\B}$ on
  $\cO_\epsilon[G][z^{-1}]/\mathfrak{m}\cO_\epsilon[G][z^{-1}],$ with Nakayama automorphism induced also from
  $\nu_{\mathcal{O}}.$

  Now suppose that $z$ is not in $\mathfrak{m}_g$. Then
  \begin{eqnarray*}
\cO_\epsilon[G][z^{-1}]/\mathfrak{m}_g\cO_\epsilon[G][z^{-1}] \cong
(\cO_\epsilon[G]/\mathfrak{m}_g\cO_\epsilon[G])[z^{-1}]
  =  \cO_\epsilon[G]/\mathfrak{m}_g\cO_\epsilon[G],
\end{eqnarray*}
using \cite[Exercise 9L]{GW} for the isomorphism, and the fact that $z$ is a unit modulo
$\mathfrak{m}_g\cO_\epsilon[G]$ for the equality. In particular, by the first paragraph of the proof,
\begin{eqnarray}
\label{conc}
\textit{the desired conclusions apply to } \cO_\epsilon[G]/\mathfrak{m}_g\cO_\epsilon[G].
\end{eqnarray}

To extend this conclusion to arbitrary $g$ in $G$ we apply the results of \cite{CL}. Recall that there is a
Poisson bracket on $\mathcal{O}[G]$, under which $G$ decomposes as a disjoint union of symplectic leaves.
Moreover, if $g,h \in G$ belong to the same symplectic leaf, then
\begin{eqnarray}
\label{symp} \cO_\epsilon[G]/\mathfrak{m}_g\cO_\epsilon[G] \cong
\cO_\epsilon[G]/\mathfrak{m}_h\cO_\epsilon[G]\end{eqnarray}
by \cite[Corollary 9.4]{CL}. In fact,
$\cO_\epsilon[G]$ is a Poisson $\cO [G]$-order in the sense of \cite[2.1]{IKPoisson}, and we can, if preferred,
quote \cite[Theorem 4.2]{IKPoisson} to obtain (\ref{symp}). By \cite[Proposition 9.3 and Proposition 8.7 (b)]{CL}
there is an action of the torus $T$ as automorphisms of $\cO_\epsilon[G],$ restricting to Poisson automorphisms of
the subalgebra $\cO[G]$ induced by right and left multiplication by $T$ on $G$, preserving the Poisson order
structure in the sense of \cite[3.8]{IKPoisson}. Therefore, if $g \in G$ and $t \in T$, then
\begin{eqnarray}
\label{wind} \cO_\epsilon[G]/\mathfrak{m}_g\cO_\epsilon[G] \cong \cO_\epsilon[G]/\mathfrak{m}_{tg}\cO_\epsilon[G]
\cong \cO_\epsilon[G]/\mathfrak{m}_{gt}\cO_\epsilon[G].
\end{eqnarray}
Since the action of $T$ preserves the leaves we can conclude from (\ref{symp}) and (\ref{wind}) that
\begin{eqnarray}
\label{comb} \cO_\epsilon[G]/\mathfrak{m}_g\cO_\epsilon[G] \cong
\cO_\epsilon[G]/\mathfrak{m}_h\cO_\epsilon[G]\end{eqnarray} if $g$ and $h$ are in the same $T$-orbit of symplectic
leaves \cite[Corollary 9.4]{CL}, \cite[4.2 and 4.3]{IKPoisson}.

Recall from \cite[Theorem A.3.2 and Theorem A.2.1]{HL} (see also \cite[Section 9.3]{CL}) that the $T$-orbits of
symplectic leaves are indexed by the elements of $W\times W$, where $W$ is the Weyl group of $G.$ To be precise,
they are the double Bruhat cells
\begin{eqnarray}\label{HL}
    X_{w_1,w_2}:=B\dot{w_1}B\cap B^-\dot{w_2}B^-,
  \end{eqnarray}
where $\dot{w_1}$, $\dot{w_2}$ are chosen from the normaliser $N_G(T)$ to represent $w_1, w_2\in W$.

Note that the localisation with respect to $z$ corresponds exactly to the localisation over the big cell $BB^-$,
as explained in \cite[proof of Theorem 7.2]{CL}. In view of (\ref{conc}) and (\ref{comb}) it is therefore enough
to show that every $T$-orbit of leaves in $G$ has non-empty intersection with the big cell. That is, by
(\ref{HL}), we must check that every double Bruhat cell $X_{w_1,w_2}$ has non-empty intersection with the big
cell. This is
 easy to verify as follows: Consider the double Bruhat cells $X_{w_1,e}=Bw_1B\cap B^-$ and $X_{e,
 w_2}=B\cap B^-w_2 B^-$. Let $a\in X_{w_1,e}$ and $b\in X_{e,
 w_2}.$ Then $ab\in B^-B\cap Bw_1B \cap  B^-w_2 B^-\subseteq B^-B\cap X_{w_1,w_2}$.
\end{proof}

\section*{Acknowledgement} We would like to thank Ami Braun for telling us about the results in
\cite{Braun} and for supplying us with a preliminary version of his paper. The third author acknowledges the
support of the EPSRC grant number GR/S14900/01. All of us benefited from the support of Leverhulme research Interchange F/00158/X (UK).
\bibliographystyle{amsplain}
\bibliography{ref7}

\providecommand{\bysame}{\leavevmode\hbox to3em{\hrulefill}\thinspace}
\providecommand{\MR}{\relax\ifhmode\unskip\space\fi MR }
\providecommand{\MRhref}[2]{%
  \href{http://www.ams.org/mathscinet-getitem?mr=#1}{#2}
}
\providecommand{\href}[2]{#2}
\begin{thebibliography}{10}

\bibitem{BF}
A.~D. Bell and R.~Farnsteiner, \emph{On the theory of {F}robenius extensions
  and its application to {L}ie superalgebras}, Trans. Amer. Math. Soc.
  \textbf{335} (1993), no.~1, 407--424.

\bibitem{Braun}
A.~Braun, \emph{On smooth {C}alabi-{Y}au algebras}, preliminary version, 2006.

\bibitem{BGoo}
K.~A. Brown and K.~R. Goodearl, \emph{Homological aspects of {N}oetherian {PI}
  {H}opf algebras and irreducible modules of maximal dimension}, J. Algebra
  \textbf{198} (1997), no.~1, 240--265.

\bibitem{BGo}
\bysame, \emph{Lectures on {A}lgebraic {Q}uantum {G}roups}, Adv. Courses in
  Math., Birkh\"auser, Basel, 2002.

\bibitem{IK}
K.~A. Brown and I.~Gordon, \emph{The ramifications of the centres: quantised
  function algebras at roots of unity}, Proc. London Math. Soc. (3) \textbf{84}
  (2002), no.~1, 147--178.

\bibitem{IKPoisson}
\bysame, \emph{Poisson orders, symplectic reflection algebras and
  representation theory}, J. Reine Angew. Math. \textbf{559} (2003), 193--216.

\bibitem{BZ}
K.~A. Brown and J.~J. Zhang, \emph{Dualizing complexes and twisted {H}ochschild
  (co)homology for noetherian {H}opf algebras}, arXiv:math.RA/0603732.

\bibitem{CGM}
S.~Caenepeel, E.~De~Groot, and G.~Militaru, \emph{{Frobenius Functors of the
  second kind}}, Comm. Algebra \textbf{30} (2002), 5359--5391.

\bibitem{ChariPresley}
V.~Chari and A.~Pressley, \emph{A {G}uide to {Q}uantum {G}roups}, Cambridge
  University Press, 1994.

\bibitem{CL}
C.~De Concini and V.~Lyubashenko, \emph{Quantum function algebra at roots of
  {$1$}}, Adv. Math. \textbf{108} (1994), no.~2, 205--262.

\bibitem{EG}
P.~Etingof and V.~Ginzburg, \emph{Symplectic reflection algebras,
  {C}alogero-{M}oser space, and deformed {H}arish-{C}handra homomorphism},
  Invent. Math. \textbf{147} (2002), no.~2, 243--348.

\bibitem{FP}
E.~M. Friedlander and B.~J. Parshall, \emph{Modular representation theory of
  {L}ie algebras}, Amer. J. Math. \textbf{110} (1988), no.~6, 1055--1093.

\bibitem{Katrin}
K.~Gehles, \emph{Properties of {C}herednik algebras and graded {H}ecke
  algebras}, Ph.D. thesis, University of Glasgow, 2006, short article version
  in preparation.

\bibitem{Gerst}
M.~Gerstenhaber and S.~D. Schack, \emph{Algebraic cohomology and deformation
  theory}, Deformation theory of algebras and structures and applications, NATO
  Adv. Sci. Inst. Ser. C, vol. 247, Kluwer Acad. Publ., 1988.

\bibitem{GW}
K.~R. Goodearl and R.~B. Warfield, Jr., \emph{An {i}ntroduction to
  {n}oncommutative {N}oetherian {r}ings}, London Mathematical Society Student
  Texts, vol.~16, Cambridge University Press, Cambridge, 1989.

\bibitem{IGG}
I.~Gordon, \emph{Baby {V}erma modules for rational {C}herednik algebras}, Bull.
  London Math. Soc. \textbf{35} (2003), no.~3, 321--336.

\bibitem{GR}
S.~Griffeth and A.~Ram, \emph{Affine {H}ecke algebras and the {S}chubert
  calculus}, European J. Combin. \textbf{25} (2004), no.~8, 1263--1283.

\bibitem{Harrison}
D.~K. Harrison, \emph{Commutative algebras and cohomology}, Trans. Amer. Math.
  Soc. \textbf{104} (1962), no.~2, 109--204.

\bibitem{HL}
T.~J. Hodges and T.~Levasseur, \emph{Primitive ideals of {${\bf C}\sb q[{\rm
  SL}(3)]$}}, Comm. Math. Phys. \textbf{156} (1993), no.~3, 581--605.

\bibitem{IR}
O.~Iyama and I.~Reiten, \emph{{F}omin-{Z}elevinsky mutation and tilting modules
  over {C}alabi-{Y}au algebras}, arXiv:math.RT/0605136.

\bibitem{Jquant}
J.~C. Jantzen, \emph{Introduction to {Q}uantum {G}roups}, Representations of
  reductive groups, Cambridge Univ. Press, 1998.

\bibitem{KadisonJones}
L.~Kadison, \emph{The {J}ones polynomial and certain separable {F}robenius
  extensions}, J. Algebra \textbf{186} (1996), no.~2, 461--475.

\bibitem{kadison}
\bysame, \emph{New examples of {F}robenius extensions}, University Lecture
  Series, vol.~14, AMS, Providence, RI, 1999.

\bibitem{Kane}
R.~Kane, \emph{Reflection {G}roups and {I}nvariant {T}heory}, CMS Books in
  Mathematics/Ouvrages de Math\'ematiques de la SMC, 5, Springer-Verlag, New
  York, 2001.

\bibitem{K}
F.~Kasch, \emph{Dualit\"atseigenschaften von {F}robenius-{E}rweiterungen},
  Math. Z. \textbf{77} (1961), 219--227.

\bibitem{KT}
H.~F. Kreimer and M.~Takeuchi, \emph{Hopf algebras and {G}alois extensions of
  an algebra}, Indiana Univ. Math. J. \textbf{30} (1981), 675--692.

\bibitem{Lusztigaffine}
G.~Lusztig, \emph{Affine {H}ecke algebras and their graded version}, J. Amer.
  Math. Soc. \textbf{2} (1989), no.~3, 599--635.

\bibitem{Luquantroot}
\bysame, \emph{Quantum groups at roots of {$1$}}, Geom. Dedicata \textbf{35}
  (1990), no.~1-3, 89--113.

\bibitem{MR}
J.C. McConnell and J.C. Robson, \emph{Noncomutative {N}oetherian {R}ings},
  Wiley-Interscience, 1987.

\bibitem{NT}
T.~Nakayama and T.~Tsuzuku, \emph{On {F}robenius extensions. {I}}, Nagoya Math.
  J. \textbf{17} (1960), 89--110.

\bibitem{Pareigis}
B.~Pareigis, \emph{Einige {B}emerkungen \"uber {F}robenius-{E}rweiterungen},
  Math. Ann. \textbf{153} (1964), 1--13.

\bibitem{Passman}
D.~S. Passman, \emph{Infinite crossed products}, Pure and Applied Mathematics,
  vol. 135, Academic Press Inc., Boston, MA, 1989.

\bibitem{RS}
A.~Ram and A.~V. Shepler, \emph{Classification of graded {H}ecke algebras for
  complex reflection groups}, Comment. Math. Helv. \textbf{78} (2003), no.~2,
  308--334.

\bibitem{Rouquier}
R.~Rouquier, \emph{Representations of rational {C}herednik algebras},
  arXiv:math.RT/0504600v2, 2005.

\bibitem{Steinberg}
R.~Steinberg, \emph{On a theorem of {P}ittie}, Topology \textbf{14} (1975),
  173--177.

\bibitem{Stcompf}
C.~Stroppel, \emph{Composition factors of quotients of the enveloping algebra
  by primitive ideals}, J. London Math. Soc. (2) \textbf{70} (2004), no.~3,
  643--658.

\bibitem{StTQFT}
\bysame, \emph{{TQFT} with corners and tilting functors in the {K}ac-{M}oody
  case}, math.RT/0605103, 2006.

\bibitem{VdB}
M.~Van~den Bergh, \emph{Existence theorems for dualizing complexes over
  noncommutative graded and filtered rings}, J. of Algebra \textbf{195} (1997),
  no.~2, 662--679.

\bibitem{Ye}
A.~Yekutieli, \emph{Dualizing complexes, {M}orita equivalence and the derived
  {P}icard group}, J. London Math. Soc. (2) \textbf{60} (1999), no.~3,
  723--746.

\bibitem{YZ}
A.~Yekutieli and J.~J. Zhang, \emph{Rings with {A}uslander dualizing
  complexes}, J. of Algebra \textbf{213} (1999), 1--51.

\end{thebibliography}
\end{document}